\documentclass[]{elsarticle}

\usepackage{lineno,hyperref}
\modulolinenumbers[5]

\journal{Journal of \LaTeX\ Templates}

\usepackage{float}
\usepackage[caption = false]{subfig}
\usepackage{amsmath}
\usepackage{booktabs}

\DeclareMathOperator{\sgn}{sgn}

\usepackage[dvipsnames]{xcolor}
\newcommand{\rv}[1]{\textcolor{black}{#1}}

\begin{document}

\begin{frontmatter}

\title{Unconditionally stable higher order \\ semi-implicit level set method for advection equations\tnoteref{mytitlenote}}
\tnotetext[mytitlenote]{The research was supported by VEGA 1/0314/23 and APVV 19-0460.}

\author[mysecondaryaddress]{Peter Frolkovi\v{c}\corref{mycorrespondingauthor}}
\cortext[mycorrespondingauthor]{Corresponding author}
\ead{peter.frolkovic@stuba.sk}

\author[mysecondaryaddress]{Nikola Gajdo\v{s}ov\'a}
\ead{nikola.gajdosova@stuba.sk}

\address[mysecondaryaddress]{Department of Mathematics and Descriptive Geometry, Slovak Technical University, Radlinsk\'eho 11, 81005 Bratislava, Slovakia}

\begin{abstract}
We present compact semi-implicit finite difference schemes on structured grids for numerical solutions of advection by an external velocity and by a speed in the normal direction that are applicable in level set methods. The recommended numerical scheme is third order accurate for the linear advection \rv{in the 2D case} with a space dependent velocity. \rv{Using a combination of analytical and numerical tools in the von Neumann stability analysis, the third order scheme is claimed to be unconditionally stable.} We also present a simple high-resolution scheme that gives a TVD (Total Variation Diminishing) approximation of the spatial derivative for the advected level set function \rv{in the 1D case}. In the case of nonlinear advection, a semi-implicit discretization is proposed to linearize the problem. The compact \rv{implicit part of the} stencil of numerical schemes contains unknowns only in the upwind direction. \rv{Consequently, algebraic solvers like the fast sweeping method can be applied efficiently to solve the resulting algebraic systems.} Numerical tests to evolve a smooth and non-smooth interface and an example with a large variation of the velocity confirm the good accuracy of the third order scheme even in the case of very large Courant numbers. \rv{The advantage of the high-resolution scheme is documented for examples where the advected level set functions contain large jumps in the gradient.}
\end{abstract}

\begin{keyword}
\texttt{level set methods \sep implicit finite differences}
\MSC[2010] 35L60 \sep65M06 
\end{keyword}

\end{frontmatter}


\section{Introduction}
\label{sec1}

Numerical methods to solve mathematical models expressed by partial differential equations are an important tool for applications of such models in research and industry. As they provide only approximate solutions, one is interested in numerical methods that are accurate and robust enough at the same time. In this paper, we attempt to offer a candidate for these types of schemes to solve a prototype of advection equations used in level set methods \cite{set99,osh02}.

The basic idea of the level set methods is to describe dynamic interfaces that can have a complex shape of evolving curves in 2D and evolving surfaces in 3D. For that purpose, a time dependent level set function is considered of which the zero level set implicitly represents in each time the position of the interface. In this way, the level set function can be determined by solving a nonlinear advection equation in which the velocity is typically prescribed by some external velocity field and/or by a speed in the normal direction. The advection equation can then be solved numerically, e.g., on a uniform structured grid using finite difference methods.

In particular, we consider the following nonlinear advection equation,
\begin{equation}
    \label{equation}
    \partial_t \phi + \left(\vec{u} +  \delta \frac{\nabla \phi}{\lvert\nabla \phi\lvert} \right) \cdot \nabla \phi = 0 \,, \quad \phi({\bf x},0)=\phi^0({\bf x}) ,
\end{equation}
where $\phi=\phi({\bf x},t)$ for ${\bf x} \in R^d$ and $t>0$ is the unknown level set function given at $t=0$ by the given function $\phi^0$. The vector field $\vec{u}=\vec{u}({\bf x})$ prescribes the movement of all level sets by an external velocity,  and $\delta$ is the speed in the normal direction given by the normalized gradient.

Numerical solutions of the level set equation \eqref{equation} are of great interest in research and applications, see monographs or review articles \cite{set99, osh02,gibou2018review} for an overview.
Concerning a large variation of the applications of level set methods, we are concerned with the tracking of interfaces in two-phase flows \cite{sussman1998improved,olsson2007conservative,frolkovic2016flux}, groundwater flow with moving water table \cite{holm_method_1999,herreros2006application, fro12,robinson2023new}, evolving porous media \cite{van2009crystal,schulz2017effective,ray2019numerical,garttner2020efficiency,kelm2022comparison}, forest fire propagation \cite{mallet2009modeling,frolkovic2015semi,alessandri2021parameter}, and image segmentation by subjective surfaces \cite{sarti2000subjective,mikula2005co,bourgine2009extraction}.

We are interested in numerical methods that do not require constraints on the choice of discretization steps to ensure the stability of computations.
Such restrictions are usually quantified by the so-called (grid) Courant numbers, which typically must be small enough to provide stable numerical results. The restriction can be unpractical to be fulfilled in cases where a large variation of Courant numbers occurs due to, e.g., the large variations of discretization steps for unfitted grids with computational domains having complex boundaries, see the so-called "small cut cells" problem in \cite{fmu15,may2017explicit,frolkovic2018semi, Engwer2020A3673, Xie2022}. Moreover, large time steps, and consequently large Courant numbers, are suitable in problems where the time dependent solution is approaching stationary form or when an auxiliary time variable is used to solve stationary problems by time marching methods or relaxation algorithms \cite{fmu15,li_absolutely_2021,hahn2022finite}.

To derive an implicit scheme with no stability restriction on time steps, we follow several techniques that are popular in proposing numerical methods for the solution of hyperbolic problems. First, we apply the so-called Lax-Wendroff (or Cauchy-Kowalevskaya) procedure in connection with finite Taylor series in time, where the time derivatives are replaced by terms involving space derivatives using the relation given by a partial differential equation. The standard form of this procedure uses terms that involve only spatial derivatives that are then approximated by some discretization methods in space \cite{qiu_finite_2003,leveque_finite_2004,toro_riemann_2009}. We follow the approach in which mixed derivatives are used \cite{zorio_approximate_2017,carrillo2019compact,carrillo2021lax} combined with the idea that the sequence of terms obtained in the Taylor series can be approximated in decreasing order of accuracy \cite{qiu_finite_2003,seal_high-order_2014,frolkovic2023high}. Using these tools, we construct the third order accurate implicit scheme that \rv{produces algebraic systems which can be solved efficiently by the solvers such as the fast sweeping method \cite{zhao2005fast}. The accuracy and stability of the scheme are studied for the linear advection equation and smooth solutions, but the method is successfully applied for the nonlinear form \eqref{equation} with nonsmooth solutions. In the latter case, the accuracy order is decreased, but the method gives better experimental orders of convergence for the chosen representative examples than when computed with some second order schemes \cite{frolkovic2018semi}. The most important property of the third order scheme is that it can be claimed to be unconditionally stable using the von Neumann stability analysis. Due to the complex form of the scheme, especially in the 2D case, this property can only be shown by a combination of analytical and numerical tools in the analysis, as previously used by other authors \cite{billett1997on,ahmed2011third,frolkovic2018semi}.}

In contrast to hyperbolic problems that describe conservation laws for which discontinuous solutions must be considered, the solutions of the nonconservative level set advection equation \eqref{equation} are supposed to be continuous. However, the gradient of the level set function can contain discontinuities, and the approximation of the gradient can play an important role in some applications of level set methods. To deal with it, we use a simple relation in the 1D case between the nonconservative advection equation for the level set function and the conservative advection equation for the spatial derivative of the level set function that is used to motivate the derivation of many numerical methods, including the one in the seminal work of Osher and Sethian \cite{osh88}. We use this property to derive a high-resolution scheme to solve \eqref{equation} based on a parametric second order scheme that can be locally limited to prevent unphysical oscillations in the approximations of gradient in the spirit of Essentially Non-Oscillatory (ENO) \cite{shu_essentially_1998,osh02} \rv{or TVD (Total Variation Diminishing) \cite{harten_class_1984,sweby1984high,kemm_comparative_2011,frolkovic2022semi} approximations}. 

In summary, we offer the semi-implicit method for the solution of \eqref{equation} based on two numerical schemes. The third order accurate scheme is given here in detail only for the two-dimensional case. It offers a very good approximation of the solution $\phi$ in \eqref{equation}, but it does not ensure a nonoscillatory approximation of $\nabla \phi$ in the case of discontinuities that could be an issue in some applications of level set methods. The high-resolution scheme is simple to implement even in multidimensional problems and offers the possibility to approximate the gradient in the spirit of ENO and related methods.

The paper is structured as follows. In Section 2 we present details of all schemes in the 1D case. In Section 3 we extend the method to several dimensions, and in Section 4 we extend it to the nonlinear case. In Section 5 we present representative numerical experiments. We conclude in Section 6 and some details on two specialized topics are given in Appendix.

\section{One-dimensional case}
\label{sec2}

For clarity of presentation, we describe the method in the one-dimensional linear case and then extend it to several dimensions and the nonlinear form in \eqref{equation}. The linear advection equation for an unknown function $\phi=\phi(x,t)$ with a given velocity function $u=u(x)$ can be written as
\begin{equation}
    \label{1dadv}
    \partial_t \phi(x,t) + u(x) \partial_x \phi(x,t) = 0 \,, \quad x \in (0,L) \,, \,\, t >0 \,.
\end{equation}
Let $x_i \in [0,L]$ (with $L$ given) and $t^n \ge 0$ be discrete spatial and temporal points with the indices $i$ and $n$ running from $0$ to given values $I$ and $N$, respectively. We restrict ourselves to a uniform spatial mesh, so $h := x_{i+1}-x_i$, with $x_0=0$ and $x_I=L$ being the boundary nodes. For simplicity, we use a uniform time step $\tau:=t^{n+1}-t^n$, but the method can be used with variable time steps. In what follows, we use the short notation $\phi_i^n:=\phi(x_i,t^n)$ and similarly for the partial derivatives of $\phi$. Analogously, $u_i:=u(x_i)$. \rv{Later, we will introduce the values of numerical solution that will be denoted by $\Phi_i^n \approx \phi_i^n$.} 

The equation \eqref{1dadv} must be accompanied by a given initial function $\phi^0=\phi^0(x)$ and given boundary functions $\phi_0=\phi_0(t)$ (if $u(0) >0$) and $\phi_L=\phi_L(t)$ (if $u(L)<0)$), which are used to define the discrete values
\begin{eqnarray*}
\label{ic}
    \phi_i^0=\phi^0(x_i) \,, \,\, i=0,1,\ldots,I \,,\\
    \label{bc}
    \phi^n_0=\phi_0(t^n) \hbox{  if  } u_0>0 \,, \,\, n=1,2,\ldots,N \,,\\
    \nonumber
    \phi^n_I=\phi_L(t^n) \hbox{  if  } u_I<0 \,, \,\, n=1,2,\ldots,N \,.
\end{eqnarray*}

We begin our study with a Taylor series expansion in a form suitable to derive an implicit type of schemes,
\begin{equation}
    \label{1dtaylor}
    \phi_i^{n-1}=\phi_i^n-\tau \partial_t \phi_i^n + \frac{\tau^2}{2} \partial_{tt} \phi_i^n - \frac{\tau^3}{6} \partial_{ttt} \phi_i^n + \mathcal{O}(\tau^4) \,.
\end{equation}
The idea of the Lax-Wendroff (or Cauchy-Kowalevskaya) procedure is to replace the time derivatives in \eqref{1dtaylor} using the equation \eqref{1dadv} by terms that contain spatial derivatives. \rv{Afterwards, finite difference approximations are used to derive a numerical scheme.} In what follows, we do it gradually. First, we derive a parametric form of the second order accurate scheme, then its high-resolution extension, and, finally, we derive the third order accurate scheme.

\subsection{Parametric second order accurate scheme}
\label{sec2a}

Unlike the original Lax-Wendroff procedure where all time derivatives in \eqref{1dtaylor} are replaced by space derivatives \cite{qiu_finite_2003,leveque_finite_2004} using \eqref{1dadv}, we use a partial Lax-Wendroff procedure where mixed derivatives are allowed \cite{zorio_approximate_2017,frolkovic2018semi,carrillo2019compact},
\begin{eqnarray}
    \label{1dlw1}
    \partial_t \phi_i^n = - u_i \partial_x \phi_i^n \,, \quad
    \partial_{tt} \phi_i^n = - u_i \partial_{tx} \phi_i^n \,.
     \end{eqnarray}
Applying \eqref{1dlw1} to \eqref{1dtaylor} we obtain
\begin{eqnarray}
    \label{1dtaylorlw2}
    \phi_i^{n-1} = \phi_i^n +  \tau u_i \partial_x \phi_i^n
     - \frac{\tau^2}{2} u_i \partial_{tx} \phi_i^n + \mathcal{O}(\tau^3) \,.
\end{eqnarray}
To obtain a fully discrete scheme, one can approximate the terms after $\tau$ and $\tau^2$ in \eqref{1dtaylorlw2} with a second and a first order accurate finite difference, respectively. 
To do so, we introduce non-dimensional Courant numbers,
\begin{equation*}
    \label{cn}
    C_i := \frac{\tau u_i}{h} \,.
\end{equation*}
First, we derive the scheme for the case $C_i > 0$ that determines the upwind form of finite difference approximations. 
Later, we present the scheme for the general case of arbitrary signs of $C_i$. 

To approximate the term in \eqref{1dtaylorlw2} after $\tau$, we consider the parametric approximation 
\begin{equation}
    \label{1dt1}
    h \partial_x \phi_i^n \approx
\phi_i^n - \phi_{i-1}^n + \frac{1-w_{i}}{2} (\phi_{i+1}^n - 2 \phi_i^n + \phi_{i-1}^n) + \frac{w_{i}}{2} (\phi_i^n - 2 \phi_{i-1}^n + \phi_{i-2}^n) \,.
\end{equation}
The approximation is second order accurate for any choice of the parameter $w_{i} \in R$, and it is third order accurate for the particular choice $w_{i}=1/3$, see \cite{wesseling2009principles, nishikawa2021truncation}. 
Note that, in general, different values of $w_i$ can be used in each time step, which we do not emphasize in the notation. 

Next, we approximate the term after $\tau^2$, where it is enough to use the first order accurate approximation. We again propose a parametric approximation, but now with the purpose of obtaining a ``compact scheme''. That is, we want to cancel the term with $\phi_{i+1}^n$ in \eqref{1dt1} for any $w_i$. To do so, we propose the following,
\begin{eqnarray*}
    \label{omega}
     \frac{\tau h}{2} \partial_{tx} \phi_i^n  \approx \frac{1}{2}\left((1-w_{i}) (\phi_{i+1}^n - \phi_{i}^n - \phi_{i+1}^{n-1} + \phi_{i}^{n-1}) \right . + \\[1ex] \nonumber \left . 
     + w_{i} (\phi_i^n - \phi_{i-1}^n - \phi_i^{n-1} + \phi_{i-1}^{n-1} ) \right) \,.
\end{eqnarray*}
Using these approximations in \eqref{1dtaylorlw2} and neglecting truncation errors when approximating $\phi_i^n$ with the values $\Phi_i^n$ of the numerical solution, we obtain the final parametric second order accurate numerical scheme for $C_i>0$,
\begin{eqnarray}
    \label{1d2oscheme}
    \Phi_i^{n} + C_i \left( \Phi_i^n - \Phi_{i-1}^n +
    \frac{1-w_{i}}{2} \left( \Phi_{i+1}^{n-1} - \Phi_{i}^{n-1} - \Phi_i^n +\Phi_{i-1}^n \right) \right. \\[1ex] \nonumber 
    \left. + \frac{w_{i}}{2} \left( \Phi_{i}^{n-1} - \Phi_{i-1}^{n-1} - \Phi_{i-1}^{n}+\Phi_{i-2}^{n} \right) \right) = \Phi_i^{n-1}
\end{eqnarray}

\rv{
If $w_{i}\neq 0$ and $w_{i}\neq 1$, the scheme \eqref{1d2oscheme} has the full stencil containing the values from $\phi_{i-2}^n$ up to $\phi_{i}^n$ in the implicit part and $\phi_{i-1}^{n-1}$ up to $\phi_{i+1}^{n-1}$ in the explicit part. The two particular choices $w_{i}=0$ and $w_{i}=1$ give the schemes with reduced stencils that can be used for approximations near boundary nodes where the full stencil is not available. For stability reasons \cite{frolkovic2018semi}, only the values $w_i\ge 0$ shall be considered, therefore, the choice $w_i>1$ in \eqref{1d2oscheme} is possible.}

The leading error term $E$ of the scheme \eqref{1d2oscheme} can be expressed in the form
\begin{eqnarray}
    \label{error2o}
    E = \frac{\tau^3}{6} \partial_{ttt} \phi_i^n+\frac{h \tau^2}{4} u_i \partial_{ttx} \phi_i^n +  \frac{h^2 \tau}{4}  (2w_i-1) u_i \partial_{txx} \phi_i^n  \\[1ex] \nonumber
   + \left. \frac{h^3}{6}  (1-3 w_i) u_i \partial_{xxx} \phi_i^n\right. \,.
\end{eqnarray}
In the case of constant velocity, that is, $u_i \equiv \bar u$, so $C_i \equiv \bar C$, one can apply the standard Lax-Wendroff procedure,
\begin{equation*}
    \label{lwbaru}
    \partial_{txx} \phi_i^n = -\bar u \partial_{xxx} \phi_i^n \,, \,\,
    \partial_{ttx} \phi_i^n = \bar u^2 \partial_{xxx} \phi_i^n \,, \,\,
    \partial_{ttt} \phi_i^n = - \bar u^3 \partial_{xxx} \phi_i^n \,.
\end{equation*}
Using it in \eqref{error2o}, we obtain
\begin{eqnarray*}
    \label{error2oBARU}
    E = \frac{h^3}{12} \bar C (1+\bar C) (2 + \bar C - 6 w_i)  \partial_{xxx} \phi_i^n \,.
\end{eqnarray*}
Clearly, the choice $w_i = (2 + \bar C)/6$ cancels the third order error term $E$, so for this choice of parameter, the scheme \eqref{1d2oscheme} is third order accurate if the velocity is constant. Such a possibility is well known also for analogous parametric fully explicit schemes \cite{wesseling2009principles,nishikawa2021truncation} or fully implicit schemes \cite{frolkovic2018semi}.

The main advantage of the compact scheme \eqref{1d2oscheme} is that the resulting linear algebraic system is defined by a lower triangular matrix, therefore, the unknowns $\Phi_i^n$ can be obtained directly if equations \eqref{1d2oscheme} are solved in the order $i=1,2,\ldots,I$. The system \eqref{1d2oscheme} must be accompanied by appropriate approximations near the boundary node $x_0$.



When deriving the scheme for $C_i<0$, we obtain
\begin{eqnarray}
    \label{1d2oschemeNeg}
    \Phi_i^{n}  
     + C_i \left( \Phi_{i+1}^n - \Phi_i^n - \frac{1-w_{i}}{2} \left( \Phi_{i+1}^{n} - \Phi_{i}^{n} - \Phi_i^{n-1} +\Phi_{i-1}^{n-1} \right) \right. \\[1ex] \nonumber
     -\left. \frac{w_{i}}{2} \left( \Phi_{i+2}^{n} - \Phi_{i+1}^{n} - \Phi_{i+1}^{n-1}+\Phi_{i}^{n-1}\right) \right)
    = \Phi_i^{n-1} \,.
\end{eqnarray}
The linear algebraic system obtained by \eqref{1d2oschemeNeg} is described by an upper triangular matrix, therefore, it can be solved directly if the equations are solved in the order $i=I-1,I-2,\ldots,0$, and if a proper treatment of approximations near the boundary node $x_I$ is used.

The general case that cover \eqref{1d2oscheme} and \eqref{1d2oschemeNeg} can be written as follows,
\begin{eqnarray}
    \label{1d2oschemegen}
    \Phi_i^{n}  
     + \lvert C_i \lvert \left( \Phi_{i}^n - \Phi_{i\mp 1}^n + \frac{1-w_{i}}{2} \left( \Phi_{i\pm 1}^{n-1} - \Phi_i^{n-1} - \Phi_{i}^{n} + \Phi_{i\mp 1}^{n}  \right) \right. \\[1ex] \nonumber
     +\left.\frac{w_{i}}{2} \left( \Phi_{i}^{n-1} - \Phi_{i\mp 1}^{n-1} - \Phi_{i\mp 1}^{n}+\Phi_{i\mp 2}^{n} \right) \right)
    = \Phi_i^{n-1} \,.
\end{eqnarray}
where $\pm=\sgn(C_i)$ and $\mp=-\sgn(C_i)$. \rv{The scheme \eqref{1d2oschemegen} was studied in \cite{frolkovic2022semi} where it was proven to be unconditionally stable using the von Neumann stability analysis for $w_i\ge 0$. Note that this analysis in \cite{frolkovic2022semi} is realized in a rigorous analytical way that we are unable to provide for the more complex numerical schemes presented later. }

To solve the linear system of algebraic equations \eqref{1d2oschemegen}, we use the fast sweeping method \cite{zhao2005fast,frolkovic2022semi} \rv{which consists of Gauss-Seidel iterations with alternating orderings of equations, the ``sweeps'', as described for \eqref{1d2oscheme} and \eqref{1d2oschemeNeg}.
In fact, the algebraic system \eqref{1d2oschemegen} can be solved exactly with two Gauss-Seidel iterations with different sweeps, if the velocity varies smoothly (e.g., linearly) between two points in the grid and changes sign only the way that $u_{\mathtt i}>0$ and $u_{\mathtt i+1}<0$ for some $\mathtt{i} \in \{1,2,\ldots,I-1\}$ . The reason is that the matrix of system \eqref{1d2oschemegen} is then reducible into two blocks, one having lower and one having upper triangular form \cite{frolkovic2022semi}.
}

\rv{
Such a property is lost if there is an index $\mathtt{i}$ such that $u_{\mathtt i}<0$ and $u_{\mathtt{i+1}}>0$, when the matrix of the system \eqref{1d2oschemegen} is not of a triangular form, and when, in general, more than two Gauss-Seidel iterations are required to solve it accurately, see more details in \cite{frolkovic2022semi}, and related numerical experiments later.
}


Concerning a choice of the value for $w_i$ in \eqref{1d2oschemegen}, we prefer the space dependent value 
\begin{equation}
    \label{wprefer}
    w_i = \frac{2+\lvert C_i\lvert}{6}
\end{equation}
that is \rv{second order accurate for variable velocity, but third order accurate in the case of constant velocity. The scheme \eqref{1d2oschemegen} then takes the form
\begin{eqnarray}
    \label{1dpreferred}
    \Phi_i^n + 
    \frac{\lvert C_i\lvert}{6}  \biggl ( 4 \Phi_i^n - 5 \Phi_{i\mp 1}^n + \Phi_{i\mp 2}^n + 2 \Phi_{i\pm 1}^{n-1} - \Phi_i^{n-1} - \Phi_{i\mp 1}^{n-1}\\[1ex] \nonumber
    + \left. \frac{\lvert C_i \lvert}{2} (\Phi_i^n - 2\Phi_{i\mp 1}^n + \Phi_{i\mp 2}^n - \Phi_{i+1}^{n-1} + 2\Phi_i^n -  \Phi_{i - 1}^{n-1}) \right) = \Phi_i^{n-1}
\end{eqnarray}
As the scheme \eqref{1dpreferred} is a special case of \eqref{1d2oschemegen}, it is unconditionally stable according to the proof in \cite{frolkovic2022semi}.
}

In the next section, we introduce a high-resolution form of the semi-implicit scheme \eqref{1d2oschemegen} where the parameter $w_i$ will depend on the numerical solution and will differ from \eqref{1dpreferred} for the grid nodes where the approximation of $\partial_x \phi_i^n$ varies significantly.


\subsection{High-resolution scheme}
\label{sec-meno}

For level set methods, the quality of the approximation for the first derivatives of the solution (the gradient) can be of great importance. The level set function itself is continuous, but the derivatives can in general be only piecewise continuous having jumps at some parts of the computational domain. Therefore, one can expect nonphysical oscillations in the approximation of the first derivative if the second order scheme \eqref{1d2oscheme} is used with a fixed stencil \cite{shu_essentially_1998,leveque_finite_2004}, that is, with a fixed value of parameters $w_i$. 

In this section, to avoid such behavior if the approximation of the gradient is important near discontinuities, we propose a nonlinear form of the scheme with parameters $w_i$ depending on the numerical solution similar to \cite{frolkovic2023high} that we adapt to the level set equation \eqref{1dadv}. Moreover, we propose the scheme in a predictor-corrector form that simplifies the solution of the resulting nonlinear algebraic equations. 

Note that the scheme \eqref{1d2oscheme} for $\Phi_i^n$ can be used to define an analogous ``conservative'' scheme for the (undivided) backward finite differences (if $C_i>0$) to approximate $\Psi_i^n \approx h \partial_x \phi(x_i,t^n)$,
\begin{equation*} \label{bwfd}
\Psi_i^n :=  \Phi_i^n - \Phi_{i-1}^n  \,, \quad i=1,2,\ldots,I \,, \\,\, n=0,1,\ldots \,\,.
\end{equation*} 
To show it, we rewrite the scheme \eqref{1d2oscheme} in the form
\begin{eqnarray}
    \label{1d2oschemepsi}
    \Phi_i^{n} + C_i \left( \Psi_i^n  +
    \frac{1-w_{i}}{2} \left( \Psi_{i+1}^{n-1} - \Psi_i^n \right) + \frac{w_{i}}{2} \left( \Psi_{i}^{n-1} - \Psi_{i-1}^{n} \right) \right) = \Phi_i^{n-1} \,.
\end{eqnarray}
Furthermore, using the notation for a ``numerical flux function'',
\begin{equation}
    \label{flux}
    F_{i} := u_i \left( \Psi_i^n +
    \frac{1}{2} \left( ({1-w_i}) \left( \Psi_{i+1}^{n-1} - \Psi_i^n \right) + w_i \left( \Psi_{i}^{n-1} -  \Psi_{i-1}^{n})\right)\right)\right)
\end{equation}
and computing the difference of \eqref{1d2oschemepsi} for $i$ and $i-1$, we obtain
\begin{eqnarray}
    \label{1d2oschemeCL}
    \Psi_i^{n} +  \frac{\tau}{h} \left( F_i - F_{i-1} \right) = \Psi_i^{n-1}  \,,
\end{eqnarray}
that can be viewed formally as a conservative finite difference scheme to solve
\begin{equation}
    \label{cl}
\partial_t \psi + \partial_x \left(u \psi\right) = 0 
\end{equation}
with $\psi:=\partial_x \phi$.
 
Compact implicit conservative schemes of the type \eqref{1d2oschemeCL} with \eqref{flux} were studied in \cite{frolkovic2023high} that we use here to define a high-resolution form of numerical fluxes in \eqref{flux} to obtain the TVD (Total Variation Diminishing) approximations of $\psi$ in \eqref{cl}. Such property in the discrete form is defined by
\begin{equation*} \label{TVD}
\sum \limits_{i=1}^I \lvert \Psi_i^n - \Psi_{i-1}^n \lvert  \le   \sum \limits_{i=1}^I \lvert \Psi_i^{n-1} - \Psi_{i-1}^{n-1}\lvert 
\end{equation*} 
if appropriate boundary conditions are supposed (e.g., periodic ones). 

An enormous amount of research is available for (TVD) high-resolution schemes in the case of fully explicit time discretizations starting with \cite{harten_class_1984,sweby1984high}, see also monographs \cite{leveque_finite_2004,toro_riemann_2009} or review in \cite{kemm_comparative_2011}. Similarly, high-resolution schemes for implicit time discretization are developed \cite{harten_class_1984,duraisamy_implicit_2007,arbogast2020third,puppo_quinpi_2022}. Here we adapt the methodology for the semi-implicit time discretization.

To propose a (nonlinear) TVD form of the parametric second order scheme \eqref{1d2oscheme} with $w_i$ depending on the numerical solution, we introduce indicators $r_i$ that measure a ratio between two variants of the second order updates in \eqref{flux} (i.e., the term multiplied by either $(1-w_i)$ or $w_i$),
\begin{equation*}
    \label{r}
    \quad r_i := \frac{\Psi_{i}^{n-1} - \Psi_{i-1}^{n}}{\Psi_{i+1}^{n-1} - \Psi_i^n} \,.
\end{equation*}
The indicators $r_i$ are clearly specific for the semi-implicit scheme and depend on the unknown values $\Psi_i^n$ (i.e., on $\Phi_i^n$). Next, we continue in the spirit of other (explicit or implicit) high-resolution methods. If we define the coefficients $s_i$,
\begin{equation*}
    \label{s}
    s_i = 1 - w_i + w_i r_i \,,
\end{equation*}
the fluxes $F_i$ can be then written in the form,
 \begin{eqnarray}
     \label{fluxhr1}
     F_i = u_i \left(\Psi_i^n + \frac{1}{2} s_i \left( \Psi_{i+1}^{n-1} - \Psi_{i}^{n}\right)\right) .
 \end{eqnarray}
The values $s_i$ in \eqref{fluxhr1} can be formally viewed as the ``slopes'' of the second order updates of the first order scheme. 

In what follows, we propose a high-resolution form of the general 1D scheme \eqref{1d2oschemegen} where the slopes $s_i$ are replaced by limited and predicted values $l_i$ that should not differ from $s_i$ whenever possible. We present here an algorithm to compute $l_i$, the motivation for such computations with a proof of the TVD property in the case of constant velocity is given in Appendix.
 
First, since the indicator $r_i$ depends on the unknown solution $\Phi_i^n$, we have to predict its value. We compute it with \eqref{1dpreferred} using the choice of the parameter $w_i$ in \eqref{wprefer}. Denoting the predicted solution by $\Phi_{i}^{n,p}$, we compute the predicted value $r_i^p$ of $r_i$,
\begin{equation}
    \label{rpredicted}
    r_i^p = \frac{\Phi_i^{n-1} - \Phi_{i\mp 1}^{n-1} -  \Phi_{i\mp 1}^n +  \Phi_{i\mp 2}^n}{\Phi_{i\pm 1}^{n-1} - \Phi_{i}^{n-1} - \Phi_{i}^{n,p}+\Phi_{i\mp 1}^{n}} \,.
\end{equation}
Next, we compute a preliminary value $l_i^p$ of $l_i$,
\begin{equation*}
    \label{l1}
    l_i^p = \max \{ 0 , \min \{ s_i^p , 2 \}\} \,, \quad s_i^p = 1 - w_i + w_i r_i^p \,. 
\end{equation*}
Finally, we compute the value $l_i$ by
\begin{equation}
    \label{l2}
    l_i = \max \{ 0, \min \{l_i^p ,  \left(\frac{2}{\lvert C_i \lvert} + l_{i\mp 1}\right) r_i^p \}\} \,.
\end{equation}
We note that other approaches can be used to compute the limited values $l_i$ of $s_i^p$ than in \eqref{l2}, see, e.g., \cite{kemm_comparative_2011,frolkovic2023high}.

Having the value $l_i$, the final scheme with one corrector step takes the form
\begin{eqnarray}
    \label{1dHRscheme}
    \Phi_i^{n} + \lvert C_i \lvert \left( \Phi_i^n - \Phi_{i\mp 1}^n +
    \frac{1}{2} l_i \left( \Phi_{i\pm 1}^{n-1} - \Phi_{i}^{n-1} - \Phi_i^{n,p} +\Phi_{i\mp 1}^n \right) \right) = \Phi_i^{n-1} .
\end{eqnarray}

\subsection{Third order accurate scheme}
\label{sec-3rd}

In this section, we use the property that the parametric approximation of $\partial_x \phi_i^n$ in \eqref{1dt1} is third order accurate for the particular choice $w_i=1/3$. With some additional effort, we extend the scheme \eqref{1d2oschemegen} to third order accuracy in space and time. 

To do so, we apply the Lax-Wendroff procedure to the derivatives $\partial_{ttt}\phi_i^n$ and $\partial_{ttx} \phi_i^n$ in \eqref{error2o} to replace them with mixed derivatives,
\begin{eqnarray}
    \label{1dlw3}
    \partial_{ttt} \phi_i^n = - u_i \partial_{ttx} \phi_i^n \,, \quad 
    \partial_{ttx} \phi_i^n =
     - \partial_{x}(u_i \partial_{tx} \phi_i^n) \,. 
\end{eqnarray}

Using now \eqref{1dlw3} together with $w=1/3$, the leading error term $E$ in \eqref{error2o} simplifies to the form,
\begin{eqnarray}
    \label{error2oafter}
    E = -u_i \frac{\tau^3}{12} \partial_{x} (u_i \partial_{tx} \phi_i^n) -\frac{h \tau^2}{12} u_i \partial_{txx} \phi_i^n   \,.
\end{eqnarray}
To obtain the third order accurate numerical scheme for $C_i>0$, we extend the second order scheme \eqref{1d2oscheme} by adding finite difference approximations of $E$ in \eqref{error2oafter}. Of course, we have to do it using the chosen stencil, that is,
\begin{eqnarray*}
    \label{1dt3a}
 u_i \frac{\tau^3}{12} \partial_{x} (u_i \partial_{tx} \phi_i^n) \approx \frac{C_i}{12} \left( C_i (\phi_i^n - \phi_{i-1}^n - \phi_i^{n-1} + \phi_{i-1}^{n-1}) \right. \\[1ex] \nonumber
\left . - C_{i-1} (\phi_{i-1}^n - \phi_{i-2}^n - \phi_{i-1}^{n-1} + \phi_{i-2}^{n-1}) \right)
 \end{eqnarray*}
 and
\begin{equation*}
    \label{1dt2aa}
    \frac{h \tau^2}{12} u_i \partial_{txx} \phi_i^n \approx
    \frac{C_i}{12}  \left(\phi_i^n - 2 \phi_{i-1}^n + \phi_{i-2}^n - \phi_i^{n-1} + 2 \phi_{i-1}^{n-1} - \phi_{i-2}^{n-1}\right) \,.
\end{equation*}



The final form of the scheme can be written as follows,
\begin{eqnarray*}
    \label{1dfinal}
    \Phi_i^n + 
    \frac{C_i}{12}  \left( 9 \Phi_i^n - 12 \Phi_{i-1}^n + 3 \Phi_{i-2}^n + 4 \Phi_{i+1}^{n-1} - 3 \Phi_i^{n-1} - \Phi_{i-2}^{n-1}\right.\\[1ex] \nonumber
    + \left.   C_i (\Phi_i^n-\Phi_{i-1}^n - \Phi_i^{n-1} + \Phi_{i-1}^{n-1}) \right. \\[1ex]\nonumber
    - \left. C_{i-1}(\Phi_{i-1}^n-\Phi_{i-2}^n - \Phi_{i-1}^{n-1} + \Phi_{i-2}^{n-1}) \right) =\Phi_i^{n-1} .
\end{eqnarray*}
The general case can be written in the following form,
\begin{eqnarray}
    \label{1dfinalgen}
    \Phi_i^n + 
    \frac{\lvert C_i\lvert}{12}  \left( 9 \Phi_i^n - 12 \Phi_{i\mp 1}^n + 3 \Phi_{i\mp 2}^n + 4 \Phi_{i\pm 1}^{n-1} - 3 \Phi_i^{n-1} - \Phi_{i\mp 2}^{n-1}\right.\\[1ex] \nonumber
    + \left. \lvert C_i \lvert (\Phi_i^n-\Phi_{i\mp 1}^n - \Phi_i^{n-1} + \Phi_{i\mp 1}^{n-1}) \right. \\[1ex]\nonumber
    \mp \left. C_{i\mp 1} (\Phi_{i\mp 1}^n-\Phi_{i\mp 2}^n - \Phi_{i\mp 1}^{n-1} + \Phi_{i\mp 2}^{n-1}) \right) =\Phi_i^{n-1} \,,
\end{eqnarray}
where $\pm=\sgn(C_i)$ and $\mp=-\sgn(C_i)$. 

We note that the scheme \eqref{1dfinalgen} in the case of constant velocity is different \rv{from the scheme \eqref{1dpreferred} because it has a larger stencil in the explicit part containing the value $\Phi_{i\mp 2}^{n-1}$. This difference makes it too complex to use analytical tools in the von Neumann stability analysis as in \cite{frolkovic2022semi} for \eqref{1dpreferred}. 
Therefore, we use a methodology that is applied for nontrivial schemes in \cite{billett1997on,ahmed2011third,frolkovic2018semi} by investigating the magnitude of amplification factor obtained from the von Neumann stability analysis for \eqref{1dfinal} using numerical tools. In particular, we define analytically, see \cite{frolkovic2018semi}, the amplification factor as a function of Courant number and variable $x \in (-\pi,\pi)$ in Wolfram Mathematica \cite{Mathematica} and inspect the maximal values of its magnitude using several graphical tools and numerical optimization procedures. As the magnitude never exceeds the value $1$ (the necessary condition for stability), we can claim with high confidence that the scheme \eqref{1dfinal} is unconditionally stable, as confirmed also by all numerical examples.}

\section{Advection in several dimensions}
\label{sec3}

The one-dimensional second order scheme from Section \ref{sec2a} and the high-resolution scheme from Section \ref{sec-meno} can be used in a straightforward manner in several dimensions applying them the dimension-by-dimension \cite{leveque_finite_2004}. We show it first for the linear advection equation
\begin{equation}
    \label{2dequation}
    \partial_t \phi({\bf x}, t) + \vec{v}({\bf x})  \cdot \nabla \phi({\bf x}, t) = 0 \,.
\end{equation}
with some given velocity function $\vec{v}$.

The partial Lax-Wendroff procedure takes in the two-dimensional case (i.e., ${\bf x}=(x,y)$, $\vec{v}=(u,v)$) the form
\begin{eqnarray*}
    \label{2dlw1}
    \partial_t \phi_{ij}^n = - u_{ij} \partial_x \phi_{ij}^n - v_{ij} \partial_y \phi_{ij}^n , \quad
    \partial_{tt} \phi_{ij}^n = - u_{ij} \partial_{tx} \phi_{ij}^n  - v_{ij} \partial_{ty} \phi_{ij}^n ,
     \end{eqnarray*}
where each term occurs analogously in $x$ and $y$ direction. We have extended the notation of Section \ref{sec2} as follows: $y_j = j h$, $j=0,1,\ldots,J$ (with $J$ given), $\phi_{ij}^n = \phi(x_i,y_j,t^n)$ and similarly for $u_{ij}$, $v_{ij}$ and $C_{ij}$. Moreover, we have to introduce the local Courant numbers for the second component of the velocity,
$$
D_{ij} = \frac{\tau v_{ij}}{h} \,.
$$

The unlimited version of the 2D scheme for $\Phi_{ij}^n \approx \phi_{ij}^n$ can be then written formally as follows,
\begin{eqnarray}
    \label{2d2oschemegen}
    \Phi_{ij}^{n}   
     + \lvert C_{i j} \lvert \left( \Phi_{i j}^n - \Phi_{i\mp 1 j}^n + \frac{1-w^x_{i j}}{2} \left( \Phi_{i\mp 1 j}^{n} - \ldots  \right) \right) \\[1ex]
     \nonumber 
    + \lvert D_{i j} \lvert \left( \Phi_{i j}^n - \Phi_{i j\mp 1}^n + \frac{1-w^y_{i j}}{2} \left( \Phi_{i j\mp 1}^{n} - \ldots  \right) \right) = \Phi_{i j}^{n-1} \,,
\end{eqnarray}
where the first term in large parentheses shall be completed analogously to \eqref{1d2oschemegen} and the same with the second one that should be adapted to the $y$ direction. The parameters $w_{ij}^x$ and $w_{ij}^y$ now correspond to $w_i$ in \eqref{1dt1} applied in the $x$ and $y$ direction, respectively.

It was shown in \cite{frolkovic2018semi} that the linear second order scheme \eqref{2d2oschemegen} is only conditionally stable \rv{using numerical tools in the von Neumann stability analysis}. For example, for the choice $w^x_{ij} = w^y_{ij} \equiv 0.5$, the scheme is stable up to Courant numbers $\lvert C_{ij}\lvert $ and $\lvert D_{ij}\lvert $ approximately equal to $7.396$ \cite{frolkovic2018semi}, which is a significant improvement compared to analogous explicit schemes. 
\rv{The maximal value of the magnitude of amplification factor is growing very slowly for larger Courant numbers reaching approximately the value $1.0454$ for Courant numbers equal $16$ for this choice of parameters \cite{frolkovic2018semi}.}

Furthermore, the high-resolution method \eqref{1dHRscheme} can be straightforwardly extended to several dimensions as follows,
\begin{eqnarray}
    \label{2dHRscheme}
    \Phi_{ij}^{n} + \lvert C_{ij} \lvert \left( \Phi_{ij}^n - \Phi_{i\mp 1 j}^n +
    \frac{1}{2} l^x_{ij} \left( \Phi_{i\pm 1 j}^{n-1} - \Phi_{i j}^{n-1} - \Phi_{ij}^{n,p} +\Phi_{i\mp 1 j}^n \right) \right) \\[1ex] \nonumber
    + \lvert D_{ij} \lvert \left( \Phi_{ij}^n - \Phi_{i j\mp 1}^n + \frac{1}{2} l^y_{ij} \left( \Phi_{i j\pm 1}^{n-1} - \Phi_{i j}^{n-1} - \Phi_{ij}^{n,p} +\Phi_{i j\mp 1}^n \right) \right)
    = \Phi_{ij}^{n-1} .
\end{eqnarray}
The predicted values $\Phi_{ij}^{n,p}$ and the values $l_{ij}^x$ and $l_{ij}^y$ of the limiter are obtained by a natural extension of the one dimensional case with $l_i$ in \eqref{l2}. \rv{We note that due to the nonlinear dependence of the limiters $l^x_{ij}$ and $l^y_{ij}$ on the numerical solution we could not provide the von Neumann stability analysis for the high-resolution scheme \eqref{2dHRscheme}. Nevertheless, the chosen numerical experiments confirm a stable behavior of numerical solutions even for very large Courant numbers that is not the case for the (unlimited) second order scheme \eqref{2d2oschemegen}, see numerical experiments in Section \ref{sec-nonsmooth} later. }

Finally, we extend the third order scheme from Section \ref{sec-3rd} for the two-dimensional case of \eqref{2dequation}. Interestingly enough, such a scheme will improve not only accuracy, but also stability. 
To derive the scheme, we have to extend \eqref{1dlw3} as follows
\begin{eqnarray*}
    \label{2dlw3two}
    \partial_{ttt} \phi_{ij}^n  = 
    - u_{ij} \partial_{ttx} \phi_{ij}^n 
     - v_{ij} \partial_{tty} \phi_{ij}^n 
     \end{eqnarray*}
together with
\begin{eqnarray*}
        \label{2dlw3help}
    \partial_{ttx} \phi_{ij}^n =
     - \partial_{x}(u_{ij} \partial_{tx} \phi_{ij}^n) - \partial_{x}(v_{ij} \partial_{ty} \phi_{ij}^n) \,, 
     \\[1ex] \label{2dlw3help2}\nonumber
     \partial_{tty} \phi_{ij}^n =
     - \partial_{y}(u_{ij} \partial_{tx} \phi_{ij}^n) - \partial_{y}(v_{ij} \partial_{ty} \phi_{ij}^n) \,.
\end{eqnarray*}

Consequently, the leading error term \eqref{error2oafter} in the 2D case takes the form
\begin{eqnarray}
    \label{2derror2oafter}
    E = - u_{ij} \frac{h \tau^2}{12} \partial_{txx} \phi_{ij}^n  - u_{ij} \frac{\tau^3}{12} \partial_{x} (u_{ij} \partial_{tx} \phi_{ij}^n)  
    - v_{ij} \frac{h \tau^2}{12} \partial_{tyy} \phi_{ij}^n \\[1ex] \nonumber
    - v_{ij} \frac{\tau^3}{12} \partial_{y} (v_{ij} \partial_{ty} \phi_{ij}^n) 
    -u_{ij} \frac{\tau^3}{12} \partial_{x} (v_{ij} \partial_{ty} \phi_{ij}^n)  -v_{ij} \frac{\tau^3}{12} \partial_{y} (u_{ij} \partial_{tx} \phi_{ij}^n)
    \,.
\end{eqnarray}
The first four terms in \eqref{2derror2oafter} also occur in the 1D case, see \eqref{error2oafter}, but the last two terms in \eqref{2derror2oafter} are specific to problems in several dimensions. 
Their finite difference approximation is rather straightforward, 
\begin{eqnarray}
    \label{2ddxv}
    u_{ij} \frac{\tau^3}{12} \partial_{x} (v_{ij} \partial_{ty} \phi_{ij}^n) \approx
    \sgn(D_{ij}) \frac{\lvert C_{ij}\lvert }{12} \left( D_{ij} (\phi_{ij}^{n} - \phi_{i j\mp 1}^{n} - \phi_{ij}^{n-1} + \phi_{i j\mp 1}^{n-1}) \right. \\[1ex] \nonumber
    - \left. D_{i\mp 1 j}  (\phi_{i\mp 1 j}^{n} - \phi_{i\mp 1 j\mp 1}^{n} - \phi_{i\mp 1 j}^{n-1} + \phi_{i\mp 1 j\mp 1}^{n-1})  \right) ,
\end{eqnarray}
where the sign in the first index, e.g. in $i\mp 1$, is decided from $\mp = - \sgn(C_{ij})$ and $\pm = \sgn(C_{ij})$, and analogously for the second index, e.g. in $j\mp 1$, one takes $\mp = - \sgn(D_{ij})$ and $\pm = \sgn(D_{ij})$. The second term is discretized analogously, 
\begin{eqnarray}
    \label{2ddyu}
    v_{ij} \frac{\tau^3}{12} \partial_{y} (u_{ij} \partial_{tx} \phi_{ij}^n) \approx
    \sgn(C_{ij}) \frac{\lvert D_{ij}\lvert }{12} \left( C_{ij} (\phi_{ij}^{n} - \phi_{i\mp 1 j}^{n} - \phi_{ij}^{n-1} + \phi_{i\mp 1 j}^{n-1}) \right. \\[1ex] \nonumber
    - \left. C_{i\mp 1 j} (\phi_{i j\mp 1}^{n} - \phi_{i\mp 1 j\mp 1}^{n} - \phi_{i j\mp 1}^{n-1} + \phi_{i\mp 1 j\mp 1}^{n-1})  \right) .
\end{eqnarray}

To define the complete scheme in the general case, we do it formally as follows,
\begin{eqnarray}
    \label{complete}
    \Phi_{ij}^{n} + 
    \frac{\lvert C_{ij} \lvert }{12} \left(9 \Phi_{ij}^n - 12 \Phi_{i\mp 1 j}^n + \ldots \right) 
    + \frac{\lvert D_{ij} \lvert }{12} \left(9 \Phi_{ij}^n - 12 \Phi_{i j \mp 1}^n + \ldots \right)  \\[1ex] \nonumber
    + \sgn(D_{ij}) \frac{\lvert C_{ij} \lvert}{12} \left( D_{ij} (\phi_{ij}^{n} - \phi_{i j\mp 1}^{n} + \ldots \right) 
    \\[1ex] \nonumber
    + \sgn(C_{ij}) \frac{\lvert D_{ij} \lvert}{12} \left(  C_{ij}  (\phi_{ij}^{n} - \phi_{i\mp 1 j}^{n} + \ldots \right) = \Phi_{ij}^{n-1} \,.
\end{eqnarray}
where the first term in parentheses of \eqref{complete} is completed as in \eqref{1dfinalgen}, analogously for the second term, but adapted to the variable $y$, and the third and fourth terms are completed according to \eqref{2ddxv} and \eqref{2ddyu}, respectively.

Note that the matrix for the system \eqref{complete} of linear algebraic equations has off-diagonal terms only in an upwind direction, therefore, algebraic solvers like the fast sweeping methods with four alternating directions of Gauss-Seidel iterations \cite{zhao2005fast} can be used efficiently. \rv{Namely, the first Gauss-Seidel iteration is realized in the order (the "sweep") with $i=1,\ldots,I$ (the outer loop) and $j=1,\ldots,I$ (the inner loop) with the initial guess for $\Phi_{ij}^n$ taking the values $\Phi_{ij}^{n-1}$. The consecutive iterations take the sweeps with one particular order reversed - the second one with $i=1,\ldots,I$ and $j=I-1,\ldots,0$, the third one with $i=I-1,\ldots,0$ and $j=I-1,\ldots,0$, and the fourth one with $i=I-1,\ldots,0$ and $j=1,\ldots,I$. 
The total number of Gauss-Seidel iterations depends on the complexity of the velocity field and the related characteristic curves \cite{zhao2005fast}, but it is supposed to be finite and independent of the refinement of the grid \cite{zhao2005fast,zhang2006high}. In our numerical experiments for 2D examples, for simplicity, we use a fixed number of Gauss-Seidel iterations, but we investigate the accuracy of results if fewer iterations are used for some examples.
}

We have applied the von Neumann stability analysis \rv{using numerical tools} for the system \eqref{complete} with constant (``frozen'') velocity values $\vec{v}$. \rv{The amplification factor is expressed in Mathematica software \cite{Mathematica} as a function of two Courant numbers and two variables $x,y \in (-\pi,\pi)$. We then investigate if the maximal magnitude of the factor does not exceed the value $1$ using some advanced numerical optimization procedures available in Mathematica \cite{Mathematica}. In this way, we can claim that the scheme \eqref{complete} is unconditionally stable with high confidence, as confirmed by all numerical experiments. }
 
\section{Nonlinear advection equation}

Until now, we have only considered the linear advection equation \eqref{2dequation}. To solve the nonlinear level set equation \eqref{equation}, we use a well-known approach of semi-implicit schemes \cite{frolkovic2015semi,boscarino_high_2016}, where a semi-linear form of PDEs is linearized by evaluating nonlinear coefficients with the values of the solution at the previous time. In particular, instead of the nonlinear equation \eqref{equation} for $t \in (t^n,t^{n+1})$ we solve the linear equation \eqref{2dequation} with
\begin{equation}
    \label{semilinear}
    \vec{v}(x,y) = \vec{u}(x,y) + \delta(x,y) \frac{\nabla \phi(x,y,t^{n-1})}{\lvert \nabla \phi(x,y,t^{n-1})\lvert} .
\end{equation}
Clearly, one loses accuracy in the numerical approximation as the normal direction of level sets is frozen in \eqref{semilinear} at the left point of the time interval. \rv{In general, the accuracy of the third order scheme \eqref{complete} can decrease down to the first order, nevertheless, for the examples with a smooth interface and appropriate choices of level set functions (e.g., the signed distance function), such a dramatic decrease is not observed, as already reported for different second order schemes \cite{fm07,mo10,hahn2019iterative} with an analogous linearization of the velocity as in \eqref{semilinear}, see also numerical examples later.}

To compute some approximations of the gradient $\nabla \phi(x_i,y_j,t^{n-1})$ in \eqref{semilinear}, one has to choose very carefully an upwind type of finite differences with appropriate accuracy.
To propose such an upwind finite difference, we follow the strategy in \cite{zhang2006high}. 
We suppose that $\phi(x,y,t)$ fulfills the standard sign property, that is, its zero level set represents a closed interface and $\phi<0$ within the closed region and $\phi>0$ otherwise \cite{set99,osh02}. 
Having such a property, we use the approximations
\begin{equation}
    \label{rt}
    h \partial_x \phi_{ij}^{n-1} \approx \left \{
    \begin{array}{lr}
    \Phi_{ij}^{n-1} - \Phi_{i-1 j}^{n-1,w^x} , \,\, & \Phi_{i-1 j}^{n-1,w^x} < \min \{ \Phi_{ij}^{n-1}, \Phi_{i+1 j}^{n-1,w^x} \}    \\[1ex]
    \Phi_{i+1 j}^{n-1,w^x} - \Phi_{ij}^{n-1} , \,\, & \Phi_{i+1 j}^{n-1,w^x} < \min \{ \Phi_{ij}^{n-1}, \Phi_{i-1 j}^{n-1,w^x} \} \\[1ex]
    0 & \hbox{otherwise}
    \end{array}
    \right.
\end{equation}
and analogously for $\partial_y \phi_{ij}^{n-1}$.

The values of $\Phi_{i\mp 1 j}^{n-1,w^x}$ are computed using the variable parametric form of the second order accurate approximation for $\partial_x \phi_i^{n-1}$ analogously to \eqref{1dt1} together with the idea of Weighted Essentially Non-Oscillatory (WENO) approximations as used in \cite{zhang2006high}. 
Namely, the value of $w^x=w^x_{ij}$ is computed as
\begin{equation*}
    \label{w}
    w^x_{ij} = \frac{1}{1+2 (r_{ij}^x)^2} \,,
\end{equation*}
where the indicators $r_{ij}^x$, and consequently the parameters $w^x_{ij}$, are computed differently for $\Phi_{i-1 j}^{n-1,w^x}$ and $\Phi_{i+1 j}^{n-1,w^x}$. Namely,
\begin{eqnarray*}
    \label{wenop}
    \Phi_{i\pm 1 j}^{n-1,w^x} = \Phi_{ij}^{n-1} \pm \frac{1-w^x_{ij}}{2} (\Phi_{i+1 j}^{n-1} - \Phi_{i-1 j}^{n-1})
    + \left.  \frac{w^x_{ij}}{2} (- 3 \Phi_{i j}^{n-1} + 4 \Phi_{i\pm 1 j}^{n-1} - \Phi_{i\pm 2 j}^{n-1}) \right. 
\end{eqnarray*}
 and
\begin{equation*}
    \label{rp1}
   r_{ij}^x = \frac{\epsilon + (\Phi_{i\pm 2 j}^{n-1} - 2 \Phi_{i\pm 1 j}^{n-1} + \Phi_{i j}^{n-1})^2}{\epsilon + (\Phi_{i+1 j}^{n-1} - 2 \Phi_{i j}^{n-1} + \Phi_{i-1 j}^{n-1})^2} \,.
\end{equation*}
The parameter $\epsilon$ has a small value to avoid a division by zero, e.g., $\epsilon=10^{-7}$. Analogous definitions are used to define the approximation of $\partial_y \phi_{ij}^{n-1}$ in \eqref{semilinear}.

Once we have the approximation \eqref{rt}, we can evaluate the velocity $\vec{v}$ in \eqref{semilinear} at each grid point $(x_i,y_j)$, and the schemes in Section \ref{sec3} can be applied straightforwardly.

\section{Numerical Experiments}\label{sec5}

In this section, we illustrate the properties of the proposed numerical schemes on several test problems. If an exact solution is available, we use it to set the boundary conditions and the initial condition. To check the Experimental Order of Convergence (EOC), we use the exact values of the solution not only at the boundary points with the inflow boundary conditions but also at the neighboring points outside of the computational domain if necessary. The implementation is realized in Matlab software \cite{MATLAB:2020}. 

The main purpose of experiments is to show that the schemes produce good accuracy when Courant numbers are large (i.e., significantly larger than typical restrictions of explicit schemes) and that they preserve the expected order of convergence even for very large Courant numbers with no instabilities produced. 

\rv{We emphasize in numerical experiments two recommended schemes of the proposed semi-implicit method. First, the third order accurate one, see \eqref{1dfinalgen} in the 1D case and \eqref{complete} in the 2D case, is used for all examples. Second, the high-resolution scheme defined for 1D in \eqref{1dHRscheme} and for 2D in \eqref{2dHRscheme} is applied for some examples. The latter scheme clearly produces better results for a standard test example in 1D with large jumps in the space derivative when compared with the third order scheme. Moreover, the high-resolution scheme significantly improves the accuracy of numerical solutions in the 2D case for a non-smooth interface and large Courant numbers when compared to the (unlimited) second order scheme \eqref{2d2oschemegen}. Nevertheless, it does not produce better results than the third order scheme, therefore, we do not further study it here for the examples with smooth solutions afterward. Instead, we present in the Appendix an illustrative example in 2D with an exact solution having large jumps in the gradient, for which the third order method produces a numerical solution with unphysical oscillations that are reduced for the high-resolution method.
}

\rv{
The particular discretization steps $h$ and $\tau$ for each experiment are chosen in such a way that one could recognize visually an influence of approximation errors on the accuracy of numerical solutions on a coarser mesh, and, eventually, a significant improvement of the accuracy for a finer mesh obtained by halving the discretization steps. In all experiments, we do not observe instabilities in numerical results for even larger Courant numbers, but the accuracy is then very low, especially when coarser meshes are used.
}

\subsection{Advection in the 1D case with a smooth solution}
\label{subs-varv}
In the following experiment, 
\rv{we confirm the EOC of the third order accurate scheme \eqref{1dfinalgen} even when using large Courant numbers and test the accuracy of numerical solution for different numbers of Gauss-Seidel iterations.} 
The velocity is defined as $u(x)=\sin(x)$ and the exact solution by
\begin{equation*}
    \phi(x, t)=\sin(2\arctan(\tan(\frac{x}{2}e^{-t}))) \quad \Rightarrow \quad \phi(x,0)=\sin(x) \,.
\end{equation*}
The example is computed for $x \in [-\frac{\pi}{2},\frac{7\pi}{2}]$ and $t \in [0,2]$. Note that the velocity $u$ changes sign four times in the interval.

\begin{figure}[H]
\centering
\hspace{-.6cm}
\subfloat{\includegraphics[width = 2.55in]{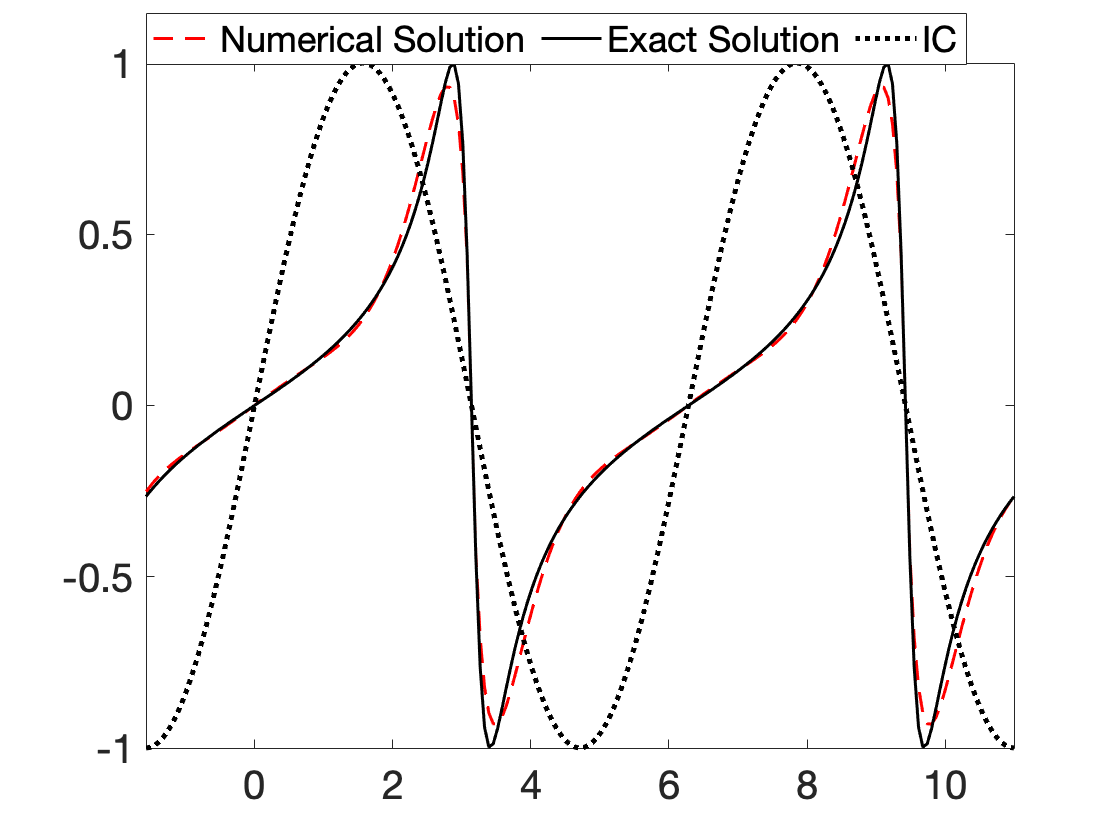}} 
\hspace{-.7cm}
\subfloat{\includegraphics[width = 2.55in]{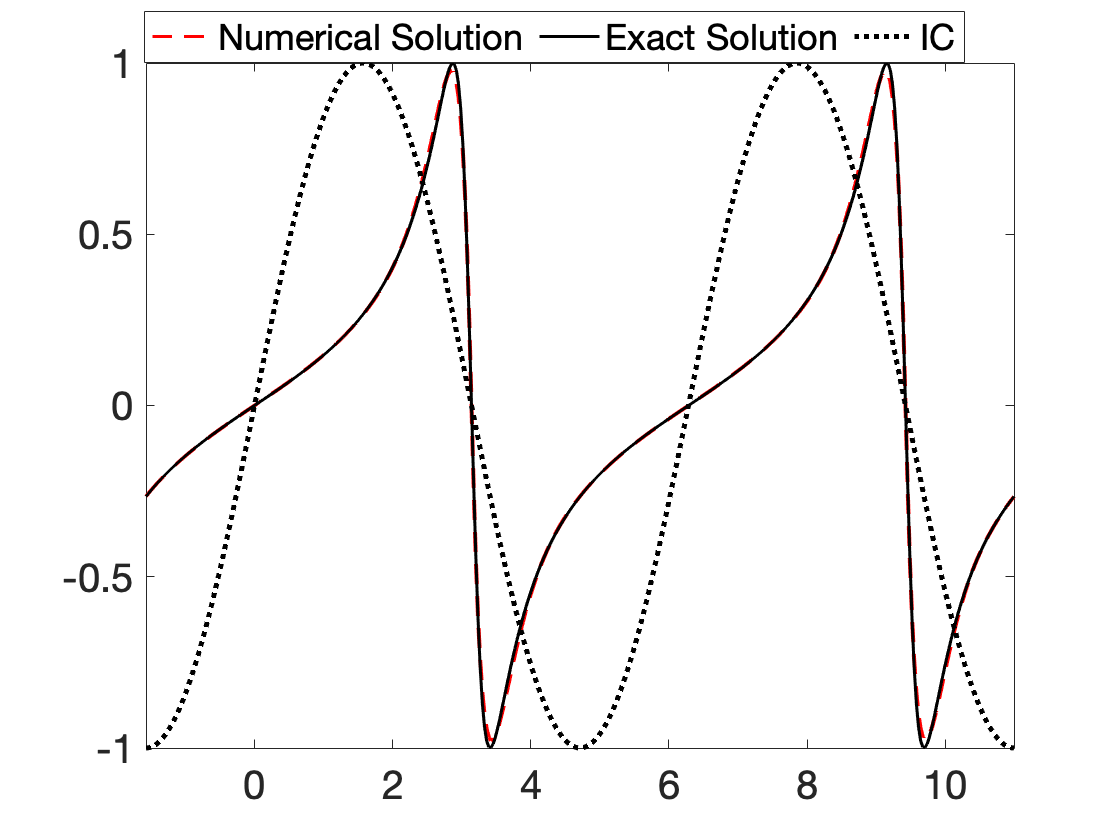}}
\caption{\rv{The example \ref{subs-varv}} - the initial function (black dotted) and the exact solution (black continuous) and the numerical solution (red dashed) at $t=2$ for $I=200$ and $N=1$ (left) and $I=400$ and $N=2$ (right) with the maximal Courant number being always approximately $32$.}
\label{1d_smooth}
\end{figure}

The error is computed by
\begin{equation}
    \label{error1}
    E_{I}^N := \tau h \sum \limits_{n=1}^{N} \sum \limits_{i=0}^I \lvert \phi_i^n - \Phi_i^n \lvert \,.
\end{equation}
 One can see in Table \ref{tab_1d_smooth_var} that the third order EOC is obtained for both the medium and large Courant numbers for sufficiently fine computational grids. \rv{Such behavior is obtained even with only two Gauss-Seidel iterations, but the accuracy is slightly improved with more iterations, especially for the coarsest mesh.}



\begin{table}
    \begin{center}
        \begin{minipage}{250pt}
        \caption{\rv{The example \ref{subs-varv}} - the errors \eqref{error1} and the EOCs for the example in Section \ref{subs-varv} for the maximal Courant number approximately $32$. \rv{The number of Gauss-Seidel iterations with different sweeps is $2$ (the third and fourth column), $4$ (the fifth and sixth column), and $6$ (the seventh and eighth column)}.}
        \label{tab_1d_smooth_var}
            \begin{tabular}{@{}llllllll@{}}
            \toprule
            $I$ & $N$ & $E_I^N$ & EOC & $E_I^N$ & EOC & $E_I^N$ & EOC  \\
            \midrule
            400 & 2 & 0.098583 & 2.70 & 0.098392 & 2.70 & 0.098340 & 2.70 \\
            800 & 4 & 0.013179 & 2.90 & 0.013162 & 2.90 & 0.013159 & 2.90 \\
            1600 & 8 & 0.001574 & 3.06 & 0.001573 & 3.06 & 0.001573 & 3.06 \\
            3200 & 16 & 0.000188 & 3.07 & 0.000188 & 3.07 & 0.000188 & 3.07 \\

            \bottomrule
            \end{tabular}
        \end{minipage}
    \end{center}    
\end{table}

\subsection{Advection in the 1D case with a nonsmooth solution}
\label{sub-nonsm}

In the following example, \rv{we compare the third order accurate scheme \eqref{1dfinalgen} with the high-resolution scheme \eqref{1dHRscheme} for a solution with discontinuous derivative.} We solve the linear advection equation \eqref{1dadv} with constant velocity $u(x)\equiv 1$ and a special form of the initial condition taken from \cite{qiushu1dspecialIC,kim2021third}, namely, $\phi(x,0)=\phi^0(x-0.5)$, where $\phi^0(x)$ is defined to be periodic and
\begin{equation*}
    \phi^0(x)= - c (x+1)+ \left \{
    \begin{array}{lr}
2\cos({\frac{3 \pi x^2}{2}})-\sqrt{3} & -1\leq x < -\frac{1}{3},
\\[1ex]
\frac{3}{2}+3\cos({2\pi x})& -\frac{1}{3}\leq x<0,
\\[1ex]
\frac{15}{2}-3\cos({2\pi x}) & 0 \leq x<\frac{1}{3},
\\[1ex]
6 \pi x(x-1) + \frac{28+4\pi+\cos({3\pi x})}{3} & \frac{1}{3}\leq x<1 ,
    \end{array} \right.
\end{equation*}
and $c=\frac{\sqrt{3}}{2}+\frac{9}{2}+\frac{2\pi}{3}$.

\begin{figure}[H]
\begin{center}
\centering
\hspace{-.7cm}
\subfloat{\includegraphics[width = 2.55in]{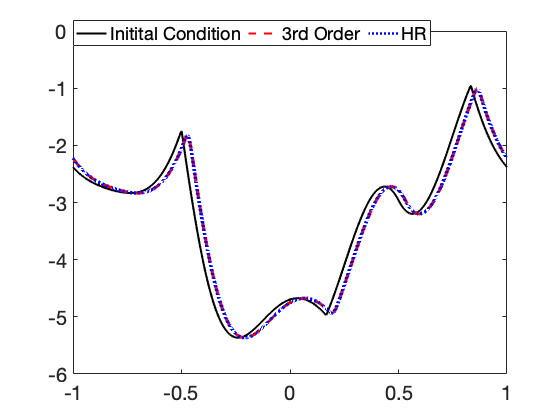}} 
\centering
\hspace{-.7cm}
\subfloat{\includegraphics[width = 2.55in]{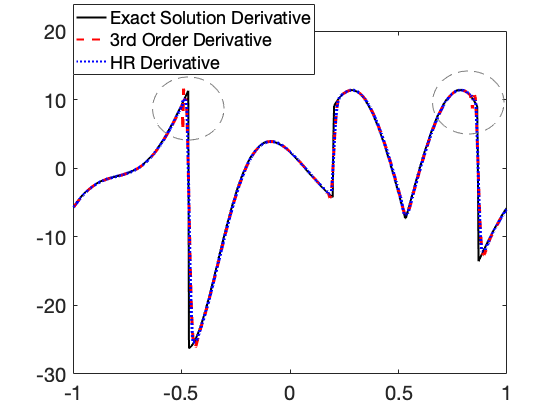}} 
\\[-2ex]
\hspace{-.7cm}
\subfloat{\includegraphics[width = 2.55in]{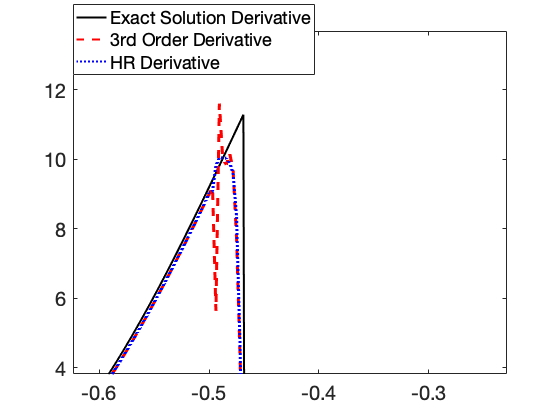}} 
\centering
\hspace{-.7cm}
\subfloat{\includegraphics[width = 2.55in]{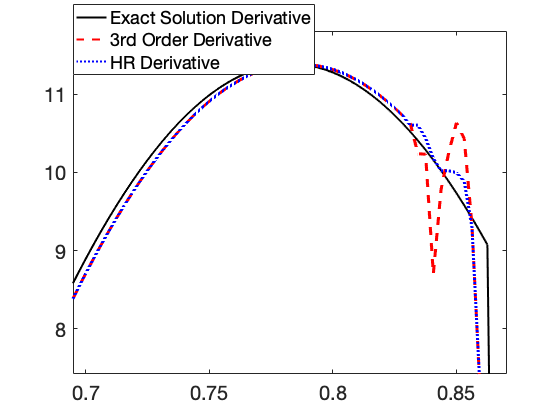}} 
\end{center}
\caption{The example \ref{sub-nonsm} - the first row contains the initial function (black full) and the numerical solutions after 2 time steps obtained with the 3rd order scheme (red dashed) and the high resolution scheme (blue dotted) for $I=640$ and $N=128$ and the Courant number equals to 5. The picture on the left contains the functions and the picture on the right contains their space derivatives approximated with the backward finite difference scheme. The second row contains the zooms of the graph in the first row.}
\label{1d_non_smooth_mp2}
\end{figure}

\begin{figure}[H]
\begin{center}
\centering
\hspace{-.7cm}
\subfloat{\includegraphics[width = 2.55in]{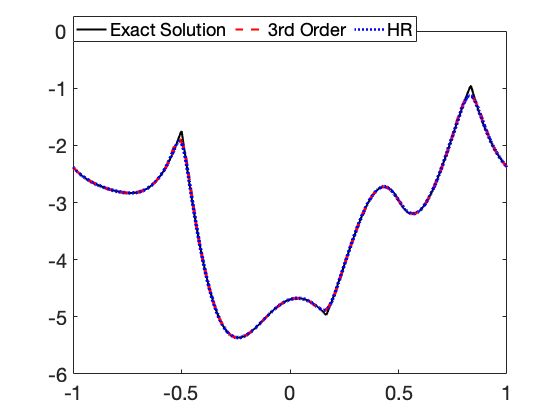}} 
\centering
\hspace{-.7cm}
\subfloat{\includegraphics[width = 2.55in]{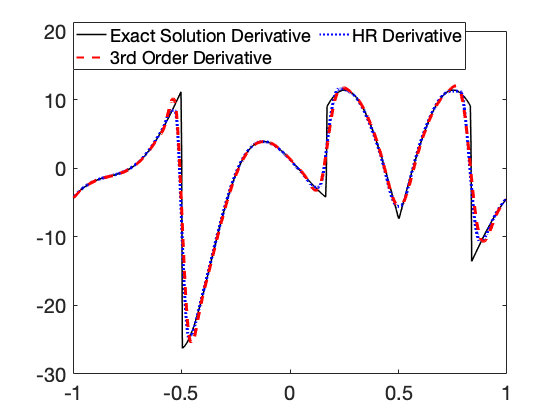}} 
\end{center}
\caption{\rv{The example \ref{sub-nonsm}} - the exact solution (black full) and the numerical solutions at $t=2$ obtained with the 3rd order scheme (red dashed) and the high resolution scheme (blue dotted) for $I=640$ and $N=128$ and the Courant number equals to 5. The left picture contains the functions, and the right one contains the space derivatives approximated with the backward finite difference scheme.}
\label{1d_non_smooth}
\end{figure}

We consider the intervals $x \in [-1,1]$ and $t \in [0,2]$. 
In Figures \ref{1d_non_smooth_mp2} and \ref{1d_non_smooth} we present the numerical solutions obtained by the high-resolution scheme \eqref{1dHRscheme} and the third order scheme \eqref{1dfinalgen}.
In Figure \ref{1d_non_smooth_mp2}, one can see that the oscillations in the approximation of $\partial_x \phi$ occur for the third order scheme at the initial time, but they are not amplified as the scheme is stable, see Figure \ref{1d_non_smooth} for the solution at $t=2$. The high-resolution method successfully reduces such oscillations at each time step. In Table \ref{tab_1d_non_smooth_var}, we present the comparison of errors \eqref{error1} and EOCs for both methods. Clearly, the high-resolution method produces smaller errors for this example. \rv{Due to the non-smooth exact solution, both schemes exhibit the EOCs below the value $2$, but still significantly larger than for a first order scheme.}

\begin{table}[H]
    \begin{center}
        \begin{minipage}{300pt}
        \caption{\rv{The example \ref{sub-nonsm}} - the errors \eqref{error1} and the EOCs for the example in Section \ref{sub-nonsm} with the Courant number equals $5$ for the 3rd order scheme (the left part) and the high-resolution scheme (the right one). }\label{tab_1d_non_smooth_var}
            \begin{tabular}{@{}lllllll@{}}
            \toprule
            $I$ & $N$ & $E_{I}$ & EOC & $E_{I}$  & EOC \\
            \midrule
            160 & 32   & 0.2900 &  1.58    & 0.2466 & 1.76 \\
            320 & 64   & 0.0874 & 1.73 & 0.0699 & 1.82 \\
            640 & 128  & 0.0281 & 1.63 & 0.0213 & 1.71 \\
            1280 & 256 & 0.0096 & 1.56 & 0.0070 & 1.6 \\
            \bottomrule
            \end{tabular}
    \end{minipage}
    \end{center}
\end{table}


\subsection{Advection in the 2D case}
\label{num-2d}

In this section, we present several linear and nonlinear test problems in the two-dimensional case of \eqref{equation}. \rv{All of them are solved using the third order accurate scheme \eqref{complete}. Additionally, the first example with a nonsmooth solution is solved with the second order scheme \eqref{2d2oschemegen} using $w^x_{ij} = w^y_{ij} = 0.5$ and the high-resolution scheme \eqref{2dHRscheme} to show a significant improvement in the accuracy for the latter scheme. Afterwards, we solve the remaining examples only with the third order scheme as it gives the most accurate results for the chosen examples. In Appendix, we present a comparison of the third order scheme and the high resolution scheme for an illustrative example for which the advected level set function is a complex distance function having large jumps in the gradient away from its zero level set.
}

To solve the resulting linear systems of algebraic equations, we use the fast sweeping method \cite{zhao2005fast} with \rv{ a fixed number iterations for all used grids as described at the end of Section \ref{sec3}. For convenience, we always use eight Gauss-Seidel iterations and we compare the accuracy of numerical solutions if only four iterations are used for some examples. 
}

\begin{figure}[H]
\centering
\hspace{-.7cm}
\subfloat{\includegraphics[width = 2.55in]{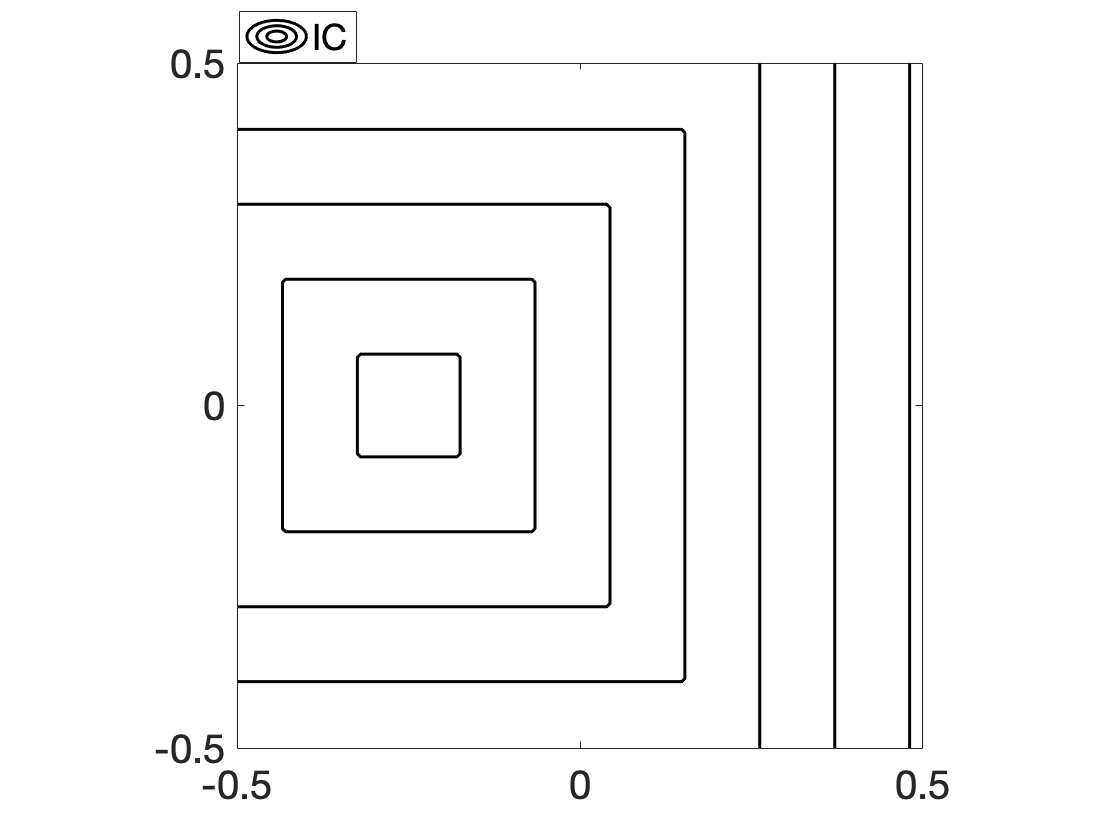}} 
\hspace{-.7cm}
\subfloat{\includegraphics[width = 2.55in]{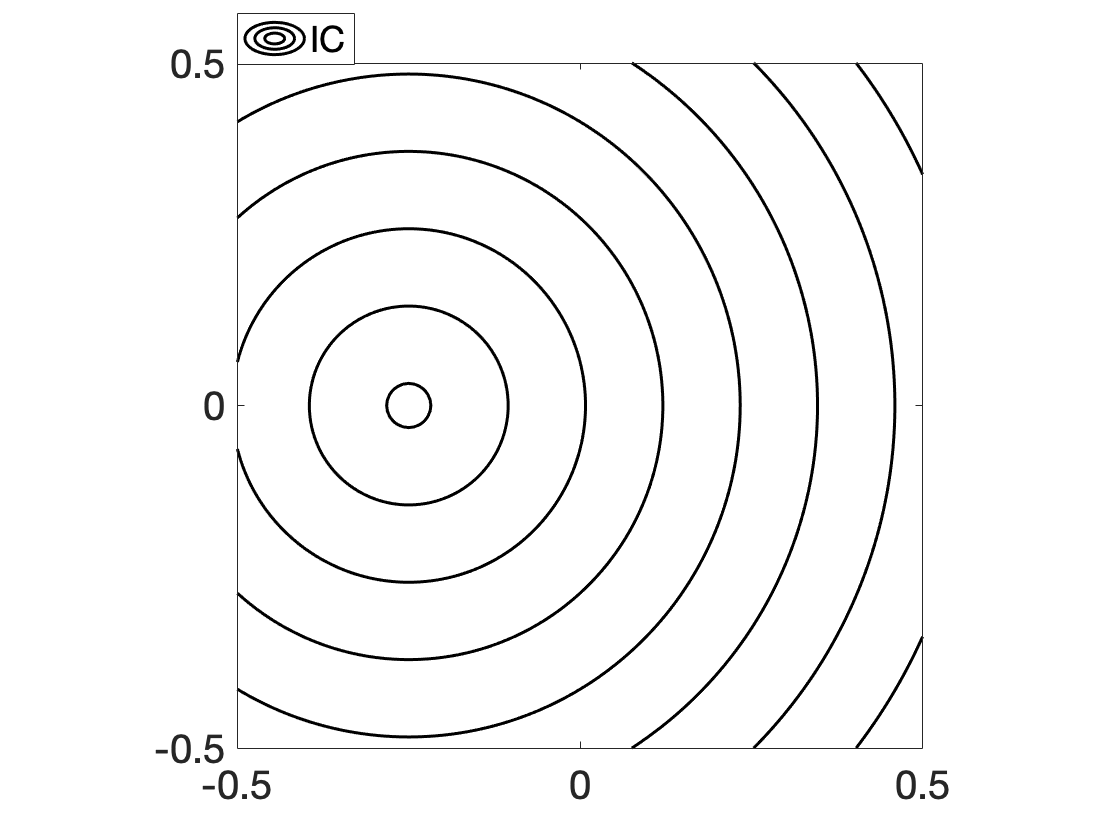}}
\caption{\rv{The examples in Section \ref{num-2d}} - the initial condition for the nonsmooth square level sets (left) and the smooth circular level sets (right).}
\label{fig-initial}
\end{figure}

Firstly, we present examples of \eqref{equation} with the velocity given as a sum of a linear velocity field $\vec{u}=(-y,x)$ and a nonlinear one describing the movement in the normal direction, namely,
\begin{equation}
\begin{pmatrix}
-y\\
x
\end{pmatrix}
+ \delta \frac{\nabla \phi}{\lvert \nabla \phi \lvert},
\label{rotShrinkExpVel}
\end{equation}
where $\delta$ is a constant. The linear part $\vec{u}$ of the velocity describes a rotation around the origin with period $2 \pi$, and the nonlinear part describes an expansion of level sets if $\delta>0$ and a shrinking if $\delta<0$. We choose two representative initial conditions defined by level set functions for a nonsmooth and a smooth interface, namely, a square interface and a circular interface, see Figure \ref{fig-initial}. The functions are defined later within the corresponding exact solutions.

Afterward, an example with exponentially varying velocity will be presented for which the property of no restriction on the time step can be used with a clear profit.

We compute the maximal Courant number $\mathcal{C}$ as given here,
\begin{equation}
    \label{maxcn}
    \mathcal{C} = \max \{ \max \limits_{i,j,n} \, \lvert C_{ij} \lvert , \,  \max \limits_{i,j,n} \, \lvert D_{ij} \lvert \}
\end{equation}
The error is computed by
\begin{equation}
    \label{error2}
    E_{I}^N := \tau h^2 \sum \limits_{n=1}^{N} \sum \limits_{i,j=0}^I \lvert \phi_{i j}^n - \Phi_{i j}^n \lvert \,.
\end{equation}

\subsubsection{Rotation of a quartic function}
\label{exquartic}

First, we \rv{check the convergence order of the third order scheme \eqref{complete}} for a quartic initial function
with the velocity $\vec{u}$ defined in \eqref{rotShrinkExpVel} and $\delta=0$. The exact solution is given by
\begin{eqnarray}
\nonumber
    \phi(\tilde{x},\tilde{y},t) = \tilde{x}^4 + \tilde{y}^4,
    \\[1ex]
    \label{xtilde}
    \tilde{x}=x\cos({t})+y\sin({t})+ 0.25,
    \quad\tilde{y}=y\cos({t})-x\sin({t}) \,,
\end{eqnarray}
where $(x,y)\in\Omega=[-1,1]\times[-1,1]$ and $t\in[0,\pi]$. We use the exact values of the solution only at the inflow part of the boundary $\partial \Omega$. In Table \ref{tab_2d_quadratic} we can see \rv{the expected results of the EOC after four Gauss-Seidel iterations that are further improved with eight Gauss-Seidel iterations. We note that twelve iterations did not bring any further changes for the errors presented in Table \ref{tab_2d_quadratic} .}

\begin{table}[H]
    \begin{center}
        \begin{minipage}{250pt}
        \caption{\rv{The example \ref{exquartic}} - the errors \eqref{error2} and the EOCs with $\mathcal{C}\approx 16$ for the results with the 3rd order scheme obtained with 4 Gauss-Seidel iterations  (the third and fourth column) and with 8 Gauss-Seidel iterations (the last two columns).}\label{tab_2d_quadratic}
            \begin{tabular}{@{}llllll@{}}
            \toprule
            $I$ & $N$ & $E_I^N$ & EOC & $E_I^N$ & EOC \\
            80 & 8 & 0.04684 & 2.90 & 0.03912 & 3.09  \\
            160 & 16 & 0.00565 & 3.05& 0.00394 & 3.31  \\
            320 & 32 & 0.00068 & 3.07 & 0.00042 & 3.23 \\
            \bottomrule
            \end{tabular}
        \end{minipage}
    \end{center}    
\end{table}

\subsubsection{Level set function for a non-smooth interface}
\label{sec-nonsmooth}

We consider the square computational domain $\Omega=[-0.5,0.5]^2$ and the time interval $t\in[0,\pi]$. In the first version of this example, we consider the velocity in \eqref{rotShrinkExpVel} with $\delta=-\frac{0.1}{\pi}$, so that the interface rotates and shrinks.  The evolved function will have level sets that are non-smooth,
namely, the exact solution is defined by
\begin{equation}
\label{stvorecinterface}
    \phi(\tilde{x},\tilde{y}, t) = 
    \left \{ \begin{array}{lr}
\tilde{y}- \delta t  & \quad \tilde{y} \geq \mid \tilde{x}\mid,
\\
-\tilde{y}- \delta t & \quad -\tilde{y} \geq \mid \tilde{x}\mid,
\\
\tilde{x}- \delta t  & \quad \tilde{x} \geq \mid \tilde{y}\mid,
\\
-\tilde{x}- \delta t  & \quad -\tilde{x} \geq \mid \tilde{y}\mid,
    \end{array} \right.
\end{equation}
where the transformed coordinates $(\tilde x, \tilde y)$ are defined in \eqref{xtilde}. 

\begin{figure}[H]
\hspace{-.7cm}
  \includegraphics[width = 2.7in]{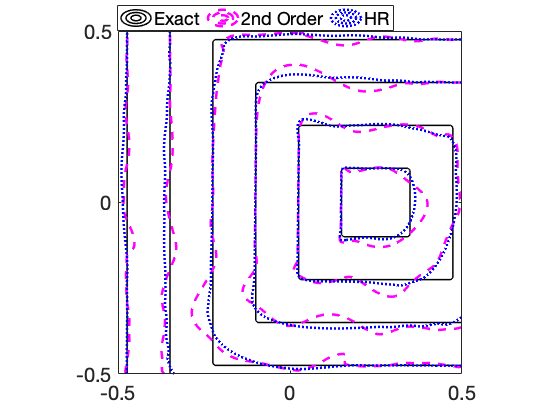}  
\hspace{-.7cm}
  \includegraphics[width = 2.7in]{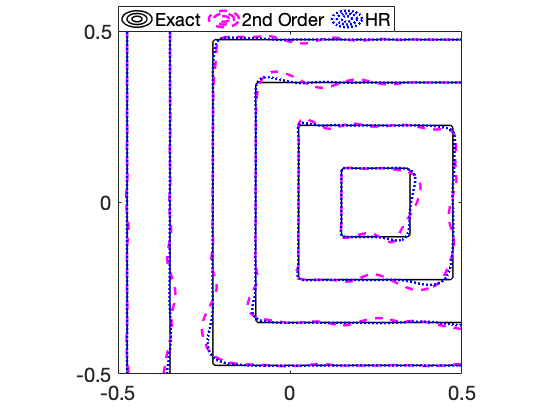}  
\caption{\rv{The example \ref{sec-nonsmooth} with the rotation and shrinking} - the exact solution (black full) and the numerical solutions obtained with the high-resolution scheme (blue dotted) and the 2nd order scheme (magenta dashed) at $t=\pi$ using $I=128$ and $N=16$ with $\mathcal{C} \approx 13.5$ (left) and $N=32$ with $\mathcal{C} \approx 6.8$ (right). \rv{The high-resolution scheme clearly reduces the oscillations produced by the second order scheme.} }
\label{2d_square_shrink_HR_2nd}
\end{figure}

\begin{figure}[H]
\hspace{-.7cm}
  \includegraphics[width = 2.7in]{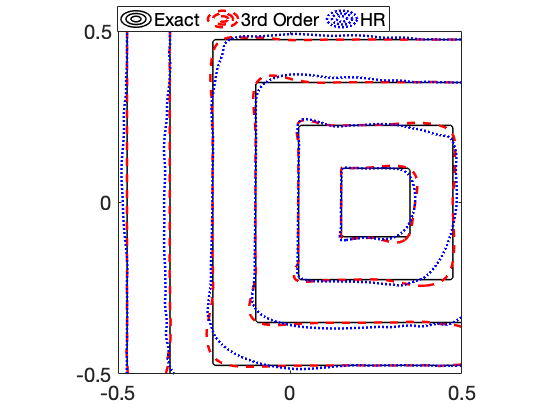}  
\hspace{-.7cm}
  \includegraphics[width = 2.7in]{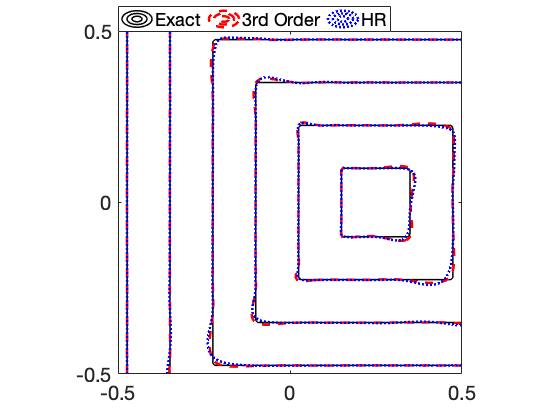}  
\caption{\rv{The example \ref{sec-nonsmooth} with the rotation and shrinking} - the exact solution (black full) and the numerical solutions obtained with the high-resolution scheme (blue dotted) and the 3rd order scheme (red dashed) at $t=\pi$ using $I=128$ and $N=16$ with $\mathcal{C} \approx 13.5$ (left) and $N=32$ with $\mathcal{C} \approx 6.8$ (right). }
\label{2d_square_shrink_HR_3rd}
\end{figure}

\rv{
The second order scheme \eqref{2d2oschemegen} with $w_{ij}^x=w_{ij}^y=0.5$ for large Courant numbers gives imprecise results that are significantly improved with the high-resolution scheme \eqref{2dHRscheme}, compare the results in Figure \ref{2d_square_shrink_HR_2nd}. One can clearly see that the oscillatory behavior of the second order scheme is significantly reduced in the results obtained by the high-resolution scheme. Furthermore, we compare the numerical results for the high-resolution scheme and the third order scheme in Figure \ref{2d_square_shrink_HR_3rd}. Clearly, the third order scheme improves the accuracy even more and it still behaves stable for large Courant numbers. 
Note that small bumps at the corners of square level sets are a typical behavior of higher order approximations \cite{fm07,mo10,saye2014high}. Finally, we compare the errors and the EOCs for the three methods in Table \ref{tab_2d_square_shrink} where the visual observation is confirmed. Note that the EOCs of all methods seem to approach approximately the first order accuracy for this example.
}

\begin{table}[H]
    \begin{center}
        \begin{minipage}{300pt}
        \caption{\rv{The first example in Section \ref{sec-nonsmooth}} - the results for the numerical solution of the rotation and shrinking of squares obtained with the 3rd order scheme (the third and fourth column), the high-resolution scheme (the fifth and sixth column) and the 2nd order scheme (the seventh and eighth column) obtained with $\mathcal{C} \approx 13.5$.}\label{tab_2d_square_shrink}
            \begin{tabular}{@{}llllllll@{}}
            \toprule
            $I$ & $N$ & $E_I^N$ & EOC & $E_I^N$ & EOC & $E_I^N$ & EOC \\
            64 & 8 & 0.019007 & 1.38 & 0.035226 & 1.12 & 0.058106 & 0.99  \\
            128 & 16 & 0.007029 & 1.44 & 0.016237 & 1.12 & 0.026536 & 1.13 \\
            256 & 32 & 0.002676 & 1.39 & 0.008180 & 0.99 &  0.011993 & 1.14  \\
            \bottomrule            
            \end{tabular}
        \end{minipage}
    \end{center}    
\end{table}


Analogously, we perform this experiment for the rotation and expansion of the initial profile by choosing $\delta=\frac{0.1}{\pi}$ in \eqref{rotShrinkExpVel}. The exact solution for this example is as follows,

\begin{equation}
\label{expandingSquareExactSol}
\begin{aligned}
    \phi(\tilde{x},\tilde{y},t) &= 
    \left \{ \begin{array}{c l}
0,& \quad   I_1
\\
-\tilde{x}-\delta t, & \quad  I_2
\\
\tilde{x}-\delta t, & \quad  I_3
\\
-\tilde{y}-\delta t, & \quad   I_4
\\
\tilde{y}-\delta t, & \quad  I_5
\\
\sqrt{(\tilde{x}-d_6)^2+(\tilde{y}-d_6)^2}-\delta t+d_6 , & \quad  I_6 \setminus \left( I_3 \cup I_5 \right)
\\
\sqrt{(\tilde{x}-d_7)^2+(\tilde{y}+d_7)^2}-\delta t+d_7, & \quad  I_7 \setminus \left( I_3 \cup I_4 \right)
\\
\sqrt{(\tilde{x}+d_8)^2+(\tilde{y}+d_8)^2}-\delta t+d_8, & \quad  I_8 \setminus \left( I_2 \cup I_4 \right)
\\
\sqrt{(\tilde{x}+d_9)^2+(\tilde{y}-d_9)^2}-\delta t+d_9, & \quad  I_9 \setminus \left( I_2 \cup I_5 \right)
    \end{array} \right.
\\
\text{where}\\
I_1&=\tilde{x}^2+\tilde{y}^2\le (\delta t)^2\\
I_2&=(\tilde{x}\le-\delta t) \,\&\, (\tilde{x}+\delta t \le \tilde{y}) \,\&\, (\tilde{y} \le -\tilde{x}-\delta t)\\
I_3&=(\tilde{x}\ge \delta t) \,\&\, (-\tilde{x}+\delta t \le \tilde{y}) \,\&\, (\tilde{y} \le \tilde{x}-\delta t))\\
I_4&=(\tilde{y} \le-\delta t) \,\&\, (\tilde{y}+\delta t \le \tilde{x})\,\&\,(\tilde{x} \le -\tilde{y}-\delta t))\\
I_5&=(\tilde{y} \ge \delta t) \,\&\, (-\tilde{y}+\delta t \le \tilde{x}) \,\&\, (\tilde{x} \le \tilde{y}-\delta t)\\
I_6&=(\tilde{x} > 0) \,\&\, (\tilde{y} > 0), \,\, d_6=\frac{1}{2}\left(\tilde{x}+\tilde{y}-\sqrt{2(\delta t)^2-(\tilde x - \tilde y)^2}\right)\\
I_7&=(\tilde{x} > 0) \,\&\, (\tilde{y} < 0), \,\, d_7=\frac{1}{2}\left(\tilde{x}-\tilde{y}-\sqrt{2(\delta t)^2-(\tilde x + \tilde y)^2}\right)\\
I_8&=(\tilde{x} < 0) \,\&\, (\tilde{y} < 0), \,\, d_8=\frac{1}{2}\left(-\tilde{x}-\tilde{y}-\sqrt{2(\delta t)^2-(\tilde x - \tilde y)^2}\right)\\
I_9&=(\tilde{x} < 0) \,\&\, (\tilde{y} > 0), \,\, d_9=\frac{1}{2}\left(-\tilde{x}+\tilde{y}-\sqrt{2(\delta t)^2-(\tilde x + \tilde y)^2}\right).
 \end{aligned}
\end{equation}

\begin{figure}[H]
\hspace{-.7cm}
  \includegraphics[width=2.7in]{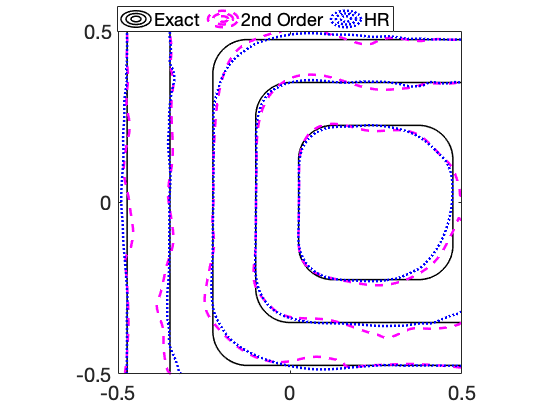}  
\hspace{-.7cm}
  \includegraphics[width=2.7in]{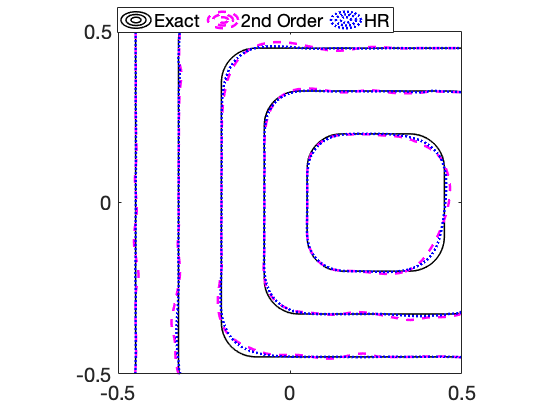}  
\caption{\rv{The example \ref{sec-nonsmooth} with the rotation and expansion} - the exact solution (black full) and the numerical solutions obtained with the high-resolution scheme (blue dotted) and the 2nd order scheme (magenta dashed) at $t=\pi$ using $I=128$ and $N=16$ with $\mathcal{C} \approx 13.5$ (left) and $N=32$ with $\mathcal{C} \approx 6.8$ (right). \rv{The high-resolution scheme clearly reduces the oscillations produced by the second order scheme.}}
\label{2d_square_expand_HR_2nd}
\end{figure}

\begin{figure}[H]
\hspace{-.7cm}
  \includegraphics[width=2.7in]{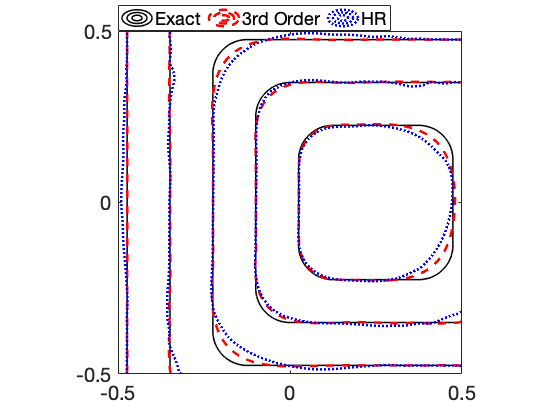}  
\hspace{-.7cm}
  \includegraphics[width=2.7in]{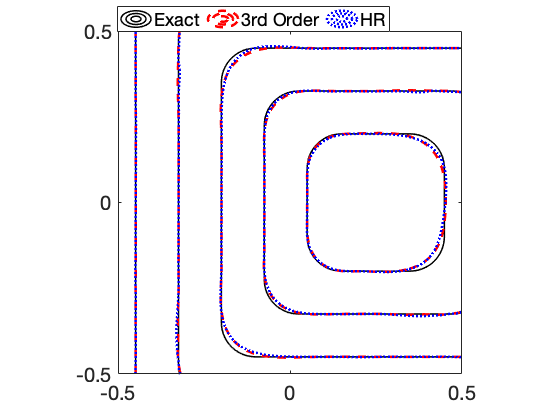}  
\caption{\rv{The example \ref{sec-nonsmooth} with the rotation and expansion} - the exact solution (black full) and the numerical solutions obtained with the high-resolution scheme (blue dotted) and the 3rd order scheme (red dashed) at $t=\pi$ using $I=128$ and $N=16$ with $\mathcal{C} \approx 13.5$ (left) and $N=32$ with $\mathcal{C} \approx 6.8$ (right). }
\label{2d_square_expand_HR_3rd}
\end{figure}

The visual comparisons are presented in Figure \ref{2d_square_expand_HR_2nd} and  \ref{2d_square_expand_HR_3rd} and the errors with the corresponding EOCs are given in Table \ref{tab_2d_square_expand}. One can again observe that the high-resolution method significantly decreases irregularities in numerical solutions obtained by the second order method and that the third order method gives improved and stable results even for large Courant numbers. \rv{The EOCs of all methods seem to approach approximately the first order accuracy for this example.}

\begin{table}[H]
    \begin{center}
        \begin{minipage}{300pt}
        \caption{\rv{The second example in Section \ref{sec-nonsmooth}} - the results for the numerical solution of the rotation and expansion of squares obtained with the 3rd order scheme (the third and fourth column), the high-resolution scheme (the fifth and sixth column) and the 2nd order scheme (the seventh and eighth column) obtained with $\mathcal{C} \approx 13.5$.}\label{tab_2d_square_expand}
            \begin{tabular}{@{}llllllll@{}}
            \toprule
            $I$ & $N$ & $E_I^N$ & EOC & $E_I^N$ & EOC & $E_I^N$ & EOC \\
            64 & 8 & 0.019387 & 1.21 &  0.035078 & 1.19 & 0.054436 & 1.02 \\
            128 & 16 & 0.008243 & 1.23 & 0.017392 & 1.01  & 0.024402 & 1.16  \\
            256 & 32 & 0.003588 & 1.20 & 0.010027 & 0.79 &  0.010599 & 1.20  \\
            \bottomrule           
            \end{tabular}
        \end{minipage}
    \end{center}    
\end{table}

\subsubsection{Level set function for a smooth interface }
\label{circle}

\rv{Next we treat analogous examples as in the previous section, but with smooth interfaces, therefore, we test only the third order accurate scheme \eqref{complete}}.
Again, we consider $\Omega=[-0.5,0.5]^2$ and $t\in[0,\pi]$. The initial condition is a distance function to the point $(-0.25,0)$,
and the exact solution is defined by
\begin{equation*}
 \begin{gathered}
    \phi(\tilde{x},\tilde{y},t)=\max \{0,\sqrt{\tilde{x}^2 + \tilde{y}^2} - \delta t \} \,,
\end{gathered}
\end{equation*}
where $\tilde x$ and $\tilde y$ are defined in \eqref{xtilde}.

In the first version, we choose $\delta=-\frac{0.1}{\pi}$, so that the initial level sets rotate and shrink.
\rv{In Table \ref{tab_2d_circle_shrink}, one can see that the third order scheme exhibits for this example the EOC approaching the value $2$ from above even for twice as large Courant numbers than for the example with non-smooth interface in the previous section. }

\begin{table}[H]
    \begin{center}
        \begin{minipage}{350pt}
        \caption{\rv{The first example in Section \ref{circle}} - the results for the numerical solution of the rotation and shrinking of circles obtained with the 3rd order scheme for $\mathcal{C}\approx 27$  (the third and fourth column) and for $\mathcal{C}\approx 13.5$ (the sixth and seventh column).}\label{tab_2d_circle_shrink}
            \begin{tabular}{@{}lllllllllll@{}}
            \toprule
            $I$ & $N$ & $E_I^N$ & EOC & $N$ & $E_I^N$ & EOC \\
            64 & 4 & 0.01023750 & 2.05 & 8 & 0.00260837 & 2.22 \\
            128 & 8 & 0.00226571 &  2.18 & 16 & 0.00056540 & 2.21  \\
            256 & 16 & 0.00052603 & 2.1 & 32 & 0.00012995 & 2.12  \\
            \bottomrule
            \end{tabular}
        \end{minipage}
    \end{center}    
\end{table}

\begin{figure}[H]
\centering
\hspace{-.7cm}
\subfloat{\includegraphics[width = 2.55in]{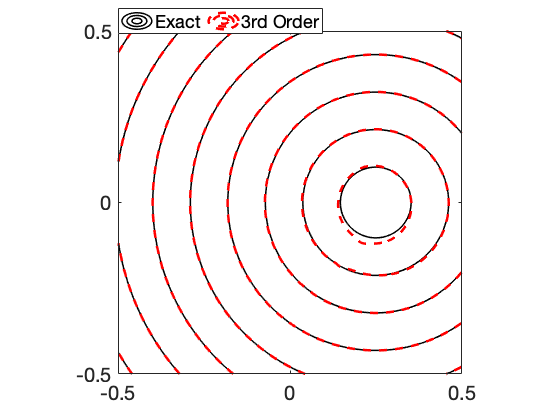}} 
\hspace{-.7cm}
\subfloat{\includegraphics[width = 2.55in]{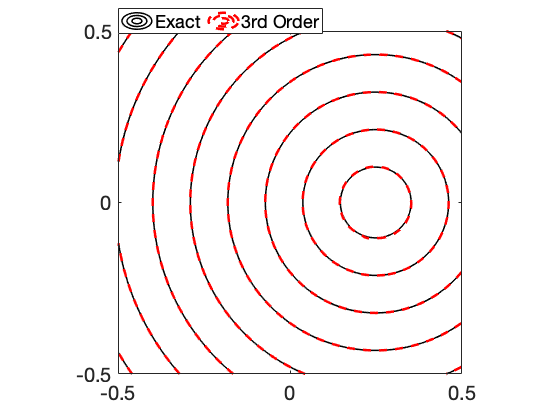}}
\caption{\rv{The example \ref{circle} with the rotation and shrinking} - the exact solution (black full) and the numerical solutions obtained with the 3rd order scheme (red dashed) at $t=\pi$ using $I=128$ and $N=8$ with $\mathcal{C} \approx 27$ (left) and $N=16$ with $\mathcal{C} \approx 13.5$ (right).}
\label{2dCircShrink}
\end{figure}

\begin{figure}[H]
\centering
\hspace{-.7cm}
\subfloat{\includegraphics[width = 2.55in]{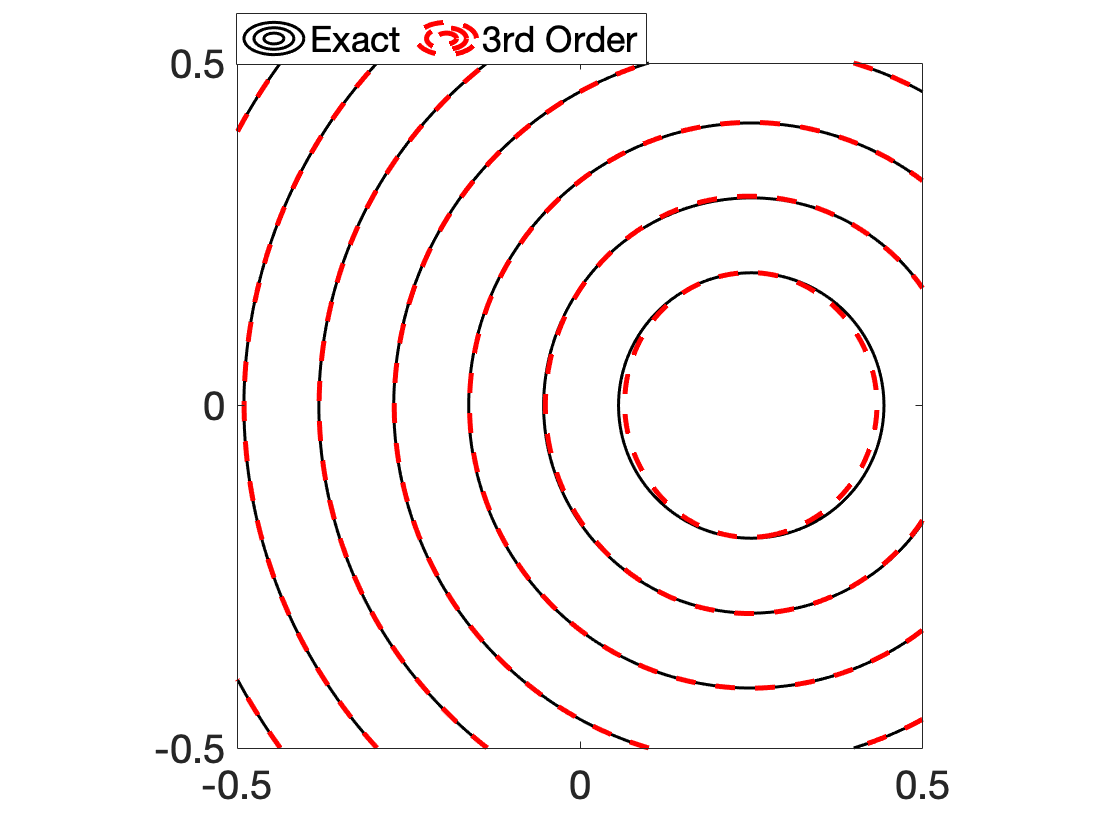}} 
\hspace{-.7cm}
\subfloat{\includegraphics[width = 2.55in]{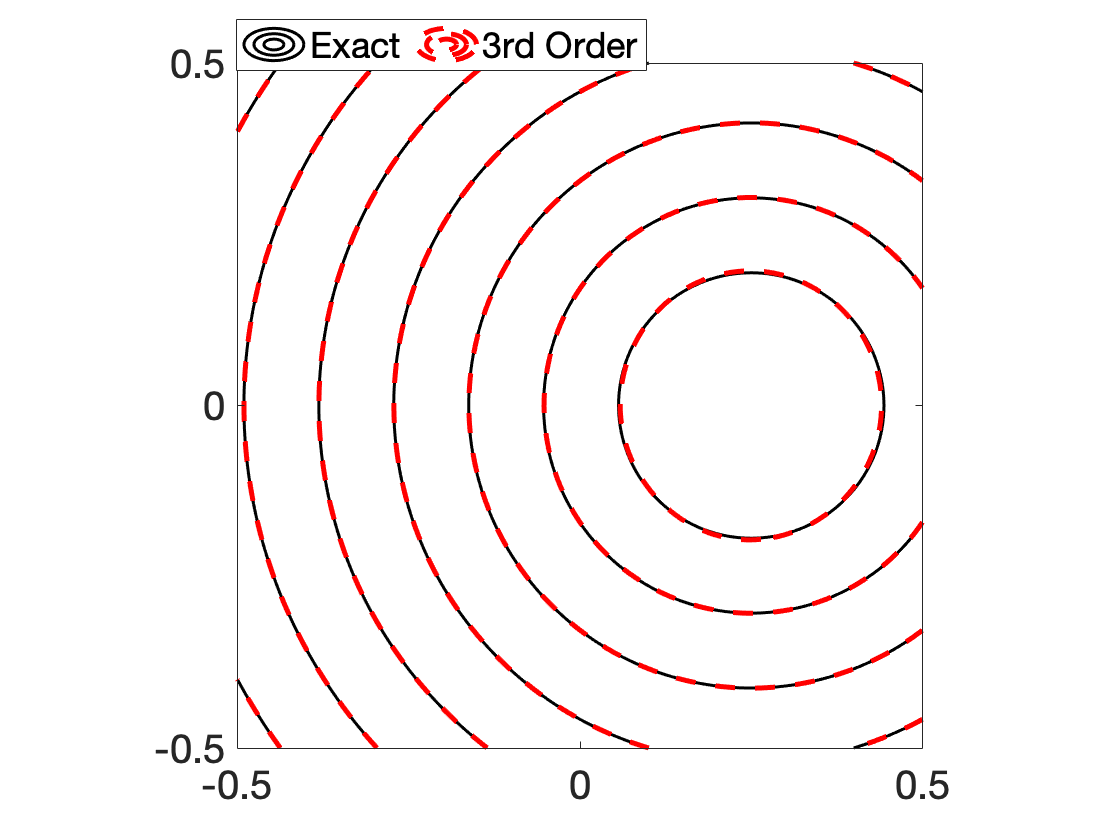}} 
\caption{\rv{The example \ref{circle} with the rotation and expansion} - the exact solution (black full) and the numerical solutions obtained with the 3rd order scheme (red dashed) at $t=\pi$ using $I=128$ and $N=8$ with $\mathcal{C} \approx 27$ (left) and $N=16$ with $\mathcal{C} \approx 13.5$ (right).}
\label{2dcircexpand}
\end{figure}

Next, we perform this experiment for the rotation and expansion of level sets by choosing $\delta=\frac{0.1}{\pi}$ in \eqref{rotShrinkExpVel}. \rv{The results are presented in Figure \ref{2dcircexpand} and in Table \ref{tab_2d_circle_expand}, where again a good accuracy of numerical solutions is confirmed even for very large Courant numbers.}

\begin{table}[H]
    \begin{center}
        \begin{minipage}{350pt}
        \caption{\rv{The second example in Section \ref{circle}} - the results for the numerical solution of the rotation and expansion of circles obtained with the 3rd order scheme for $\mathcal{C}\approx 27$  (the third and fourth column) and for $\mathcal{C}\approx 13.5$ (the sixth and seventh column).
        }\label{tab_2d_circle_expand}
            \begin{tabular}{@{}lllllllllll@{}}
            \toprule
            $I$ & $N$ & $E_I^N$ & EOC & $N$ & $E_I^N$ & EOC \\
            64 & 4 & 0.01465536 & 1.58 & 8 & 0.00578463 & 1.67 \\
            128 & 8 & 0.00467344 &  1.65 & 16 & 0.00177062 & 1.71 \\
            256 & 16 & 0.00146051 & 1.68 & 32 & 0.00054617 & 1.70  \\
            \bottomrule
            \end{tabular}
        \end{minipage}
    \end{center}    
\end{table}

\subsubsection{Exponentially varying  velocity}
\label{sec-exp}

In the following example, we illustrate the behavior of the \rv{third order} scheme \eqref{complete} for a solution of the advection equation when the velocity significantly changes its values in the computational domain. Namely, we choose the velocity that varies exponentially, 
\begin{equation*}
    \vec{v}=(u(x,y), u(x,y))=(e^{2(y-x)}, e^{2(y-x)}).
\end{equation*}
The velocity $\vec{v}$ is constant along each diagonal given by $y-x=c$ with any constant $c \in R$.

\begin{figure}[H]
\centering
\includegraphics[width = 2.55in]{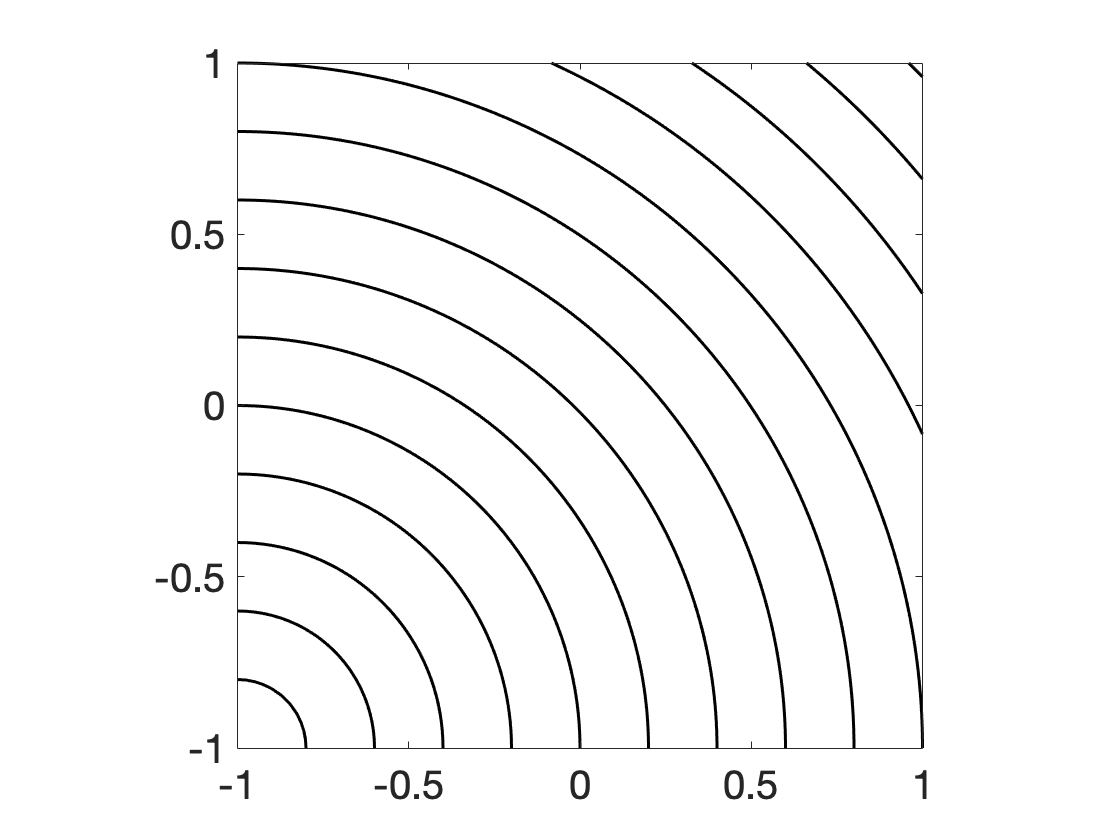}
\centering
\caption{The initial condition for the example with exponentially varying velocity.}
\label{2dexponIC}
\end{figure}

The example is defined for the square domain $\Omega=[-1,1]^2$ and the time interval $t\in[0,0.4]$. The initial condition $\phi^0$, see Figure \ref{2dexponIC}, is the distance function to the point $[-1,-1]$ defined as
\begin{equation}
\label{ICexp}
    \phi^0(x,y) = \sqrt{(x+1)^2+(y+1)^2}.
\end{equation}

\begin{figure}[H]
\centering
\hspace{-.7cm}
\subfloat{\includegraphics[width = 2.55in]{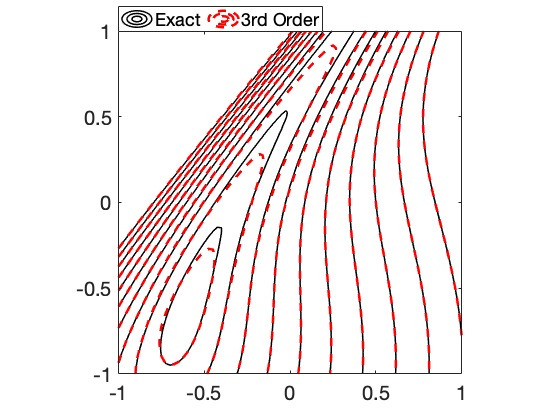}} 
\hspace{-.7cm}
\subfloat{\includegraphics[width = 2.55in]{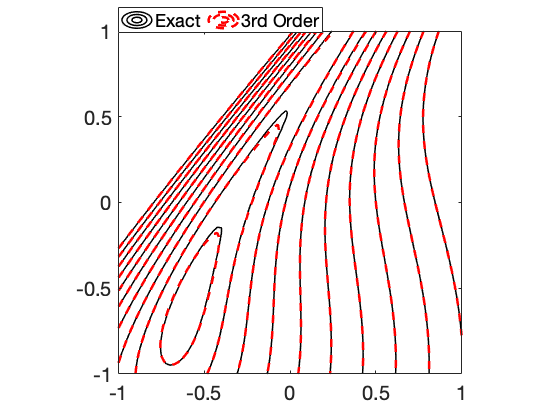}} 
\caption{\rv{The example \ref{sec-exp} with time dependent Dirichlet boundary conditions} - the results with the exact solution (black) given by \eqref{exactexp} and the numerical solutions obtained with the 3rd order scheme (red dashed) at $t=0.4$ using $I=80$ and $N=2$ with $\mathcal{C} \approx 436$  (left) and $N=8$ with $\mathcal{C} \approx 109$ (right).}
\label{2dexponinflowdbc}
\end{figure}

The exact solution $\phi=\phi(x,y,t)$ defined for $(x,y) \in R^2$ and $t\in R$  is given by
\begin{equation}
\label{exactexp}
    \phi(x,y,t) = \phi^0(x-t u(x,y),y-t u(x,y)).
\end{equation}
In the first version of the example, we set time dependent Dirichlet boundary conditions with the values given by $\phi$ from \eqref{exactexp} only at the inflow edges of the square domain. The comparison of the exact solution at the final time $t=0.4$ with numerical solutions is presented in Figure \ref{2dexponinflowdbc}, the norms of errors and the EOCs are given in Table \ref{tab_2d_exponential_inflow_bc}. One can see stable results \rv{even for very large Courant numbers and an appropriate behavior of EOCs}. Note that in this version of the example, the solution is smooth, except in one point where it takes the value $0$.

\begin{table}[H]
    \begin{center}
        \begin{minipage}{350pt}
        \caption{\rv{The first example in Section \ref{sec-exp} with time dependent Dirichlet boundary conditions} - the results for the numerical solutions obtained with the 3rd order scheme for $\mathcal{C}\approx 10.9$ (the third and fourth column) and for $\mathcal{C}\approx 109$ (the sixth and seventh column) and for $\mathcal{C}\approx 436$ (the ninth and tenth column).}\label{tab_2d_exponential_inflow_bc}
            \begin{tabular}{@{}llllllllll@{}}
            \toprule
            $I$ & $N$ & $E_I^N$ & EOC & $N$ & $E_I^N$ & EOC & $N$ & $E_I^N$ & EOC  \\
            80 & 80 & 0.0002724 & 2.57 & 8 & 0.002476 & 2.30 & 2 & 0.01849 &  2.06 \\
            160 & 160 & 0.0000435 & 2.65 & 16 & 0.0004282 & 2.53 & 4 & 0.003556 & 2.38   \\
            320 & 320 & 0.0000071 & 2.61 & 32 & 0.0000648 & 2.72 & 8 & 0.000614 & 2.54 \\

            \bottomrule
            
            \end{tabular}
        \end{minipage}
    \end{center}    
\end{table}

\begin{figure}[H]
\centering
\hspace{-.7cm}
\subfloat{\includegraphics[width = 2.55in]{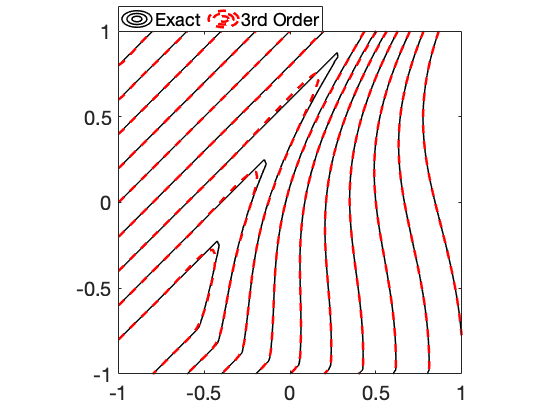}} 
\hspace{-.7cm}
\subfloat{\includegraphics[width = 2.55in]{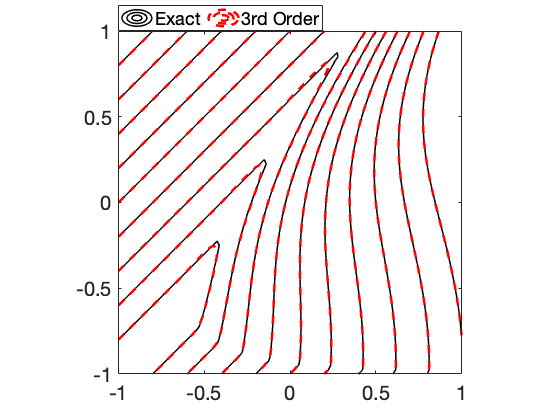}} 
\caption{\rv{The example \ref{sec-exp} with time independent Dirichlet boundary conditions} - the results with the exact solution (black) given by \eqref{exactstationary} and the numerical solutions obtained with the 3rd order scheme (red dashed) at $t=0.4$ using $I=80$ and $N=8$ with $\mathcal{C} \approx 109$  (left) and $N=80$ with $\mathcal{C} \approx 10.9$ (right).}
\label{2dexponfixedbc}
\end{figure}

In the second version of the example, we use time independent Dirichlet boundary conditions at the left and bottom edges of the domain with the values given by $\phi^0$. The solution $\phi$ is given by
\begin{equation}
\label{exactstationary}
    \phi(x,y,t) = 
    \left \{ \begin{array}{llcr}
    \begin{aligned}   
\phi^0(x-t u(x,y),y-t u(x,y), & \quad  y\ge x \text{ \& } x-t u(x,y) \ge -1
\\
\phi^0(-1,y-x-1), & \quad  y\ge x \text{ \& } x-t u(x,y) < -1
\\
\phi^0(x-t u(x,y),y-t u(x,y), & \quad  x\ge y \text{ \& } y-t u(x,y) \ge -1
\\
\phi^0(x-1-y,-1), & \quad  x\ge y \text{ \& } y-t u(x,y) < -1
\end{aligned}
    \end{array} \right.
\end{equation}
The solution reaches a stationary form in a finite time with the stationary values equilibrated much faster in the part of the square domain above its diagonal where the Courant numbers are large, see Figure \ref{2dexponfixedbc}. The third order scheme can compute the results
with a similar precision for $\mathcal{C} \approx 10.9$ and $\mathcal{C} \approx 109$. The norms of errors and the corresponding EOCs are presented in Table \ref{tab_2d_exponential_fixed_bc}. Note that the exact solution is non-smooth in this case.

\begin{table}[H]
    \begin{center}
        \begin{minipage}{350pt}
        \caption{\rv{The second example in Section \ref{sec-exp} with time independent Dirichlet boundary conditions} - the results for the numerical solutions obtained with the 3rd order scheme for $\mathcal{C}\approx 10.9$ (the third and fourth column) and for $\mathcal{C}\approx 109$ (the sixth and seventh column).}\label{tab_2d_exponential_fixed_bc}
            \begin{tabular}{@{}lllllll@{}}
            \toprule
            $I$ & $N$ & $E_I^N$ & EOC & $N$ & $E_I^N$ & EOC  \\
            80 & 80 & 0.001219 & 1.55 & 8 & 0.004463 & 1.56 \\
            160 & 160 & 0.000407 & 1.52 & 16 & 0.001508 &  1.56 \\
            320 & 320 & 0.000145 & 1.53 & 32 & 0.000521 & 1.56 \\
            \bottomrule
            
            \end{tabular}
        \end{minipage}
    \end{center}    
\end{table}

\section{Conclusions}

We present the semi-implicit method for the numerical solution of the level set equation described by the advection equation with the velocity defined by an external velocity field and by a speed in the normal direction. We propose two numerical schemes to be used with the method that have the fully upwind stencils in their implicit parts, so specialized solvers such as the fast sweeping method can be used efficiently to solve the resulting algebraic systems.

The first scheme is based on a limiting of the parametric second order scheme and can be applied straightforwardly dimension-by-dimension to problems in several dimensions. This high-resolution scheme provides in the 1D case with constant velocity the approximation of the space derivative that has provably the Total Variation Diminishing (TVD) property. The scheme has a nonlinear form due to the dependence of its parameter on the numerical solution, but the resulting algebraic system can be linearized using the predictor-corrector approach. The semi-implicit method with this scheme is suitable to numerically solve the advection equation for level set functions with large jumps in the gradient.

The second scheme is shown to be third order accurate for the linear advection equation with space dependent velocity. The scheme is unconditionally stable with high confidence in the sense of von Neumann stability analysis, which is realized using numerical tools instead of analytical ones. The scheme has a more complex form than the second order scheme and we present it here in the two-dimensional case. For chosen numerical examples where the level set function describes implicitly the advection of piecewise smooth interfaces, the third order scheme gives very good results even for very large Courant numbers violating significantly the CFL conditions of typical explicit schemes. In the future, we plan to extend the third order scheme to some high-resolution form.

\section{Appendix}

In the following two sections, we give more details for specific topics on the high-resolution scheme.

\subsection{TVD property for advection with constant velocity in 1D}

 Here, we prove that the high-resolution scheme \eqref{1dHRscheme} is TVD for the approximation of $\partial_x \phi$ in the case of the advection equation with constant velocity. 
 
 Let $C>0$, the case with a constant negative Courant number is treated analogously. 
 First, the flux $F_i$ in \eqref{fluxhr1} can be rewritten in the form
 \begin{eqnarray*}
     F_i = u_i \left(\Psi_i^n + \frac{1}{2} s_i \left( \Psi_{i+1}^{n-1} - \Psi_{i}^{n}\right)\right) =
      u_i \left(\Psi_i^n + \frac{1}{2} \frac{s_i}{r_i} \left( \Psi_{i}^{n-1} - \Psi_{i-1}^{n}\right)\right) .
 \end{eqnarray*}
Consequently, the scheme \eqref{1d2oschemeCL} can be written in the form
 \begin{eqnarray}
     \label{atvdscheme}
     \Psi_i^n - \Psi_i^{n-1} + C \left(\Psi_i^n - \Psi_{i-1}^n \right. \\[1ex] \nonumber
     + \left. \frac{1}{2} \left(\frac{s_i}{r_i} - s_{i-1}\right)
     \left( \Psi_{i}^{n-1} - \Psi_{i-1}^{n}\right) \right) = 0 \,.
 \end{eqnarray}
Moreover, as $\Psi_{i}^{n-1} - \Psi_{i-1}^{n} = \Psi_{i}^{n} - \Psi_{i-1}^{n} - (\Psi_{i}^{n} - \Psi_{i}^{n-1})$, we can rewrite \eqref{atvdscheme} as follows,
\begin{eqnarray}
     \label{atvdschemeafter}
     \Psi_i^n - \Psi_i^{n-1} +  C \frac{1+ \frac{1}{2} \left( \frac{s_i}{r_i} - s_{i-1}\right)}{1 - \frac{C}{2} \left( \frac{s_i}{r_i} - s_{i-1}\right) } \left(\Psi_i^n - \Psi_{i-1}^n\right)
     = 0 \,.
 \end{eqnarray}

Similar schemes are studied in \cite{duraisamy_implicit_2007,puppo_quinpi_2022, frolkovic2023high} with the straightforward conclusion that the scheme \eqref{atvdschemeafter} is TVD if the coefficient before $\left(\Psi_i^n - \Psi_{i-1}^n\right)$ is nonnegative. Clearly, the coefficients $s_i$ do not ensure such property in general, therefore, we replace $s_i$ in \eqref{atvdschemeafter} with limited values $l_i$,
\begin{eqnarray}
     \label{atvdschemeafterlim}
     \left(\Psi_i^n - \Psi_i^{n-1}\right) +  C \frac{1+ \frac{1}{2} \left( \frac{l_i}{r_i} - l_{i-1}\right)}{1 - \frac{C}{2} \left( \frac{l_i}{r_i} - l_{i-1}\right) } \left(\Psi_i^n - \Psi_{i-1}^n\right)
     = 0 \,.
 \end{eqnarray}
 The required property for the coefficient is obtained if the limited values fulfill the following inequalities for any $r \in \mathcal{R}$,
\begin{eqnarray}
    \label{alimiterinequalities1}
    0 \le l_{i-1} \le 2 \,, \\[1ex]
     \label{alimiterinequalities2}
     0 \le \frac{l_i}{r} \le \frac{2}{C} + l_{i-1} \,.
 \end{eqnarray}
The definition of values $l_i$ in \eqref{l2} ensures the validity of such inequalities, consequently, the coefficient in \eqref{atvdschemeafterlim} is nonnegative and the scheme \eqref{atvdschemeafterlim} is TVD for the approximation of $\partial_x \phi$.

\subsection{Example of the advection of a complex distance function in 2D}

\begin{figure}[H]
\hspace{-.7cm}
  \includegraphics[width = 2.7in]{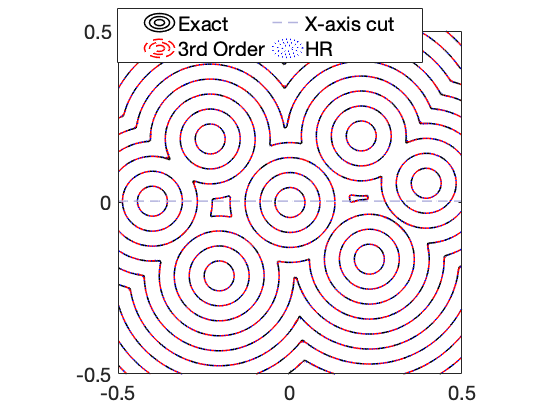}  
\hspace{-.7cm}
  \includegraphics[width = 2.7in]{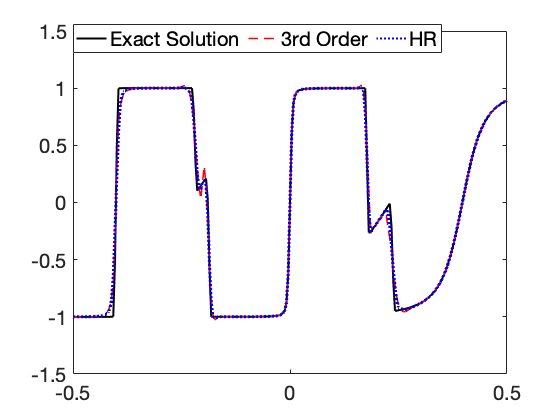}  
\caption{\rv{The example of the distance function to seven circles} - the exact solution (black full) and the numerical solutions with the high-resolution scheme (blue dotted) and the 3rd order scheme (red dashed) obtained at the second time step of the advection by the rotation where $t=0.0654$ with $I=248$ and $N=8$ with $\mathcal{C} \approx 4.1$.  The left picture contains the functions, and the right one contains the $x$-derivatives approximated with the central finite difference scheme in the cut plotted in the left picture.}
\label{2d_multiple_points_rotation}
\end{figure}

Without the purpose of going into a detailed study, we present an example for an illustration of the case for which the high-resolution scheme can bring benefits with respect to the third order scheme. We are inspired by the test example in the 1D case in Section \ref{sub-nonsm} with the solution having large jumps in the derivative. A similar situation can occur in the 2D case when the level set function takes the form of a complex distance function, as illustrated in Figure \ref{2d_multiple_points_rotation}. This function can be seen as the distance to seven circles with the smallest radius in Figure \ref{2d_multiple_points_rotation}. In the neighborhood of such an interface, the level set function is smooth, but it has large jumps in the gradient away from the interface.

We let the initial function rotate with the velocity defined in \eqref{rotShrinkExpVel} for only two time steps. The accuracy of the obtained results by the high-resolution and the third order scheme is comparable, but, analogously to Section \ref{sub-nonsm}, some oscillations occur in the approximation of the gradient with the third order scheme that are significantly reduced for the high-resolution scheme, see Figure \ref{2d_multiple_points_rotation}. Such oscillations are clearly visible for the approximation of $\partial_x \phi$ that we plot in Figure \ref{2d_multiple_points_rotation}.

\bibliography{main}

\end{document}